\documentclass[12pt,leqno]{amsart}
\usepackage[dvips]{graphicx}
\usepackage{amsfonts}
\usepackage{amssymb}
\usepackage{amsmath}
\usepackage{amscd}
\usepackage{graphicx}
\usepackage[notcite,notref]{showkeys}
\usepackage[utf8]{inputenc}   
\usepackage[T1]{fontenc}      
\usepackage{geometry}         
\usepackage[francais]{babel}  

\def\DATE{\today}
\newtheorem{theorem}{Theorem}
\newtheorem{definition}[theorem]{Definition}

\newtheorem{lemma}[theorem]{Lemma}

\newtheorem{proposition}[theorem]{Proposition}

\newcommand\R{\mathbb{R}}

\newcommand\ra{\rightarrow}

\newcommand\m{\mathfrak{m}}

\newcommand\K{\mathbb{K}}
\newcommand\Z{\mathbb{Z}}

\newcommand\p{\mathcal{P}}

\newcommand\pf{\noindent{\it Proof. }}

\newcommand\lr{\left\{ \begin{array}{l}}
\newcommand{\im}{{\rm Im}}
\def\ds{\displaystyle}
\newcommand\bu{\bullet}

\title{Depolarization and distributive laws}
\author{Elisabeth Remm}
\date{}
\address{ Universit\'{e} de Haute Alsace, LMIA 4 rue des Fr\`{e}res Lumi\`{e}%
re, 68093 Mulhouse, B.M: Universit\'{e} de Haute Alsace, IRIMAS 18 rue des Fr\`{e}res Lumi\`{e}%
re, 68093 Mulhouse }
\email{elisabeth.remm@uha.fr}

\begin{document}
\begin{abstract}
Given a vector space with two multiplications, one commutative the other anticommutative, possibly connected by a distributive law,  the depolarization principle allows to look at this triplet through a single nonassociative multiplication. This is the case of Poisson algebras. We are interested here in the cases of transposed Poisson algebras and we show in this case that depolarization cannot be done with a single multiplication. We also examine the depolarization for Hom-Lie algebras.

\end{abstract}
\maketitle
\tableofcontents

\section{Introduction and definitions}
Let $\K$ be a field of characteristic $0$ and let $V$ be a $\K$-vector space. A structure of nonassociative algebra on $V$ is given by a multiplication
$$\mu : V \otimes V \ra V$$
that we will also note $\mu(x,y)=xy$ and which satisfies axioms given by linear combinations of terms of the form $\mu(\mu(x_{\sigma(1)},x_{\sigma(2)})x_{\sigma(3)})$ and/or $\mu(x_{\sigma(1)},\mu(x_{\sigma(2)},x_{\sigma(3)})$ where $\sigma \in \Sigma_3$ is a permutation. All classical examples of algebras, such as associative, commutative, Lie and, quite surprisingly, Poisson algebras, are of this type.

To such a multiplication we can associate two other multiplications on $V$, one commutative that we will note $x \bu y$, and the other anticommutative that we will note $[x,y]$, as follows
$$
\begin{array}{l}
      x \bu y = xy+yx, \\
   \left[ x,y\right]=xy-yx. 
\end{array}
$$
So, $\mu(x,y)=xy=\frac{1}{2}(x\bu y +[x,y])$. 

The correspondence that at multiplication
$\mu$ matches the two multiplications $\bu$ and $[-,-]$ is called the polarization of $\mu$. The depolarization process is the inverse process. It associates the couple $(\bu,[-,-])$  of a commutative law and an anticommutative law, a nonassociative multiplication defined by $\mu(x,y)=xy=\frac{1}{2}(x\bu y +[x,y])$. This process can be interesting if the multiplications $\bu$ and $[-,-]$  are connected by a distributive law, that is to say a relation of the type

$\alpha_1 x_1 \bu [x_2,x_3]+\alpha_2 x_2 \bu [x_3,x_1] +\alpha_3 x_3 \bu [x_1,x_2]+\beta_1[x_1 \bu x_2,x_3]+\beta_2[x_2\bu x_3,x_1]+\beta_3 [x_3 \bu x_1,x_2]=0.
$

In this case, the study of the triple corresponding to laws $\bu$, $[-,-]$ and distributive law can be reduced to the study of a single identity relating to depolarization. This can be useful in studies of deformations, classifications or sometimes to highlight some properties of the algebra defined by such a triplet. In \cite{GR,MRPoisson}, The Poisson algebras have been studied using this point of view. Let us note also, that this approach by depolarization has shown that the operad of Lie admissible algebras is  a Koszul operad. 

In this work we will take the general framework of depolarization and focus on a class of algebras, introduced by \cite{Bai}, which is constructed as Poisson algebras from a distributive law between a Lie bracket and an associative and commutative multiplication. These algebras are called transposed Poisson algebras. First we recall the elementary principles of depolarization by explicitly describing examples. Then we resume, in a systematic way, the study of Poisson algebras. In the last section, we focus on transposed Poisson algebras and show that unlike Poisson algebras, depolarization cannot be written with a single nonassociativity rule.

\section{Distributive laws and depolarization}

The vector spaces, algebras studied in this work are always built on a  field $\K$ of charateristic zero and algebraically closed.

\subsection{Some notations}
Let $\mu$ be a multiplication on $V$ which satisfies a quadratic relation
\begin{equation}\label{mu}
\begin{array}{ll}a_1(x_1x_2)x_3+a_2(x_2x_1)x_3+a_3(x_3x_2)x_1+a_4(x_1x_3)x_2+a_5(x_2x_3)x_1+a_6(x_3x_1)x_2 & \\
+b_1x_1(x_2x_3)+b_2x_2(x_1x_3)+b_3x_3(x_2x_1)+b_4x_1(x_3x_2)+b_5x_2(x_3x_1)+b_6x_3(x_1x_2)&=0.
\end{array}
\end{equation}
for any $x_1,x_2,x_3 \in V$.  To simplify the language, instead of saying that the multiplication $\mu$ verifies a relation like the previous one, we will say also that $\mu$ verifies the axiom described by this relation.

Each permutation $\sigma \in \Sigma_3$ defines a isomorphism, noted also $\sigma$:
$$\sigma : V^{\otimes ^3} \ra V^{\otimes ^3}$$
by $\sigma(x_1,x_2,x_3)=(x_{\sigma(1)},x_{\sigma(2)},x_{\sigma(3)}$ which can be extended by
$$\sigma((x_1x_2)x_3)=(x_{\sigma(1)}x_{\sigma(2)})x_{\sigma(3)}, \ \ \sigma(x_1(x_2x_3))=x_{\sigma(1)}(x_{\sigma(2)}x_{\sigma(3)}).$$
Let $\K[\Sigma_3]$ be the algebra group associated with the symmetric group $\Sigma_3$, that is the algebra generated by $\Sigma_3$. Let $\mathcal{B}=\{Id,\tau_{1,2},\tau_{1,3},\tau_{2,3},c,c^2\}$ be a fixed ordered basis of this algebra, where $\tau_{i,j}$ is the transposition between $i$ and $j$ and $c$ the cycle $(1,2,3)$. In all this work the decompositions of the vectors of $ \K[\Sigma_3]$ will be established with respect to this base and we will no longer specify it. If $v=r_1Id+r_2\tau_{1,2}+r_3\tau_{1,3}+r_4\tau_{2,3}+r_5c+r_6c^2\in \K[\Sigma_3]$, the previous isomorphism associated with a permutation extends to the elements of $\K[\Sigma_3]$:
$$v((x_1x_2)x_3)=(r_1Id+r_2\tau_{1,2}+r_3\tau_{1,3}+r_4\tau_{2,3}+r_5c+r_6c^2)((x_1x_2)x_3)$$
the same for the right parenthesis. Using these isomorphisms, (\ref{mu})  is more simply written:
\begin{equation}\label{muv}
\begin{array}{l}v_1((x_1x_2)x_3)+v_2(x_1(x_2x_3))=0
\end{array}
\end{equation}
with $v_1=(a_1,a_2,a_3,a_4,a_5,a_6)$ and $v_2=(b_1,b_2,b_3,b_4,b_5,b_6)$ corresponding to the decomposition of $v_1$ and $v_2$ is the fixed basis $\mathcal{B}$. Conversely, the data of two vectors $(v_1=(a_i),v_2=(b_i))$ defines the relation (\ref{mu}).

Assume that $\mu$ satisfies a relation as (\ref{mu}). 
For any permutation $\sigma \in \Sigma_3$ the following relation is also true
$$
\begin{array}{ll}a_1(x_{\sigma(1)}x_{\sigma(2)})x_{\sigma(3)}+a_2(x_{\sigma(2)}x_{\sigma(1)})x_{\sigma(3)}+a_3(x_{\sigma(3)}x_{\sigma(2)})x_{\sigma(1)}+a_4(x_{\sigma(1)}x_{\sigma(3)})x_{\sigma(2)} & \\
+a_5(x_{\sigma(2)}x_{\sigma(3)})x_{\sigma(1)}+a_6(x_{\sigma(3)}x_{\sigma(1)})x_{\sigma(2)} 
+b_1x_{\sigma(1)}(x_{\sigma(2)}x_{\sigma(3)})+b_2x_{\sigma(2)}(x_{\sigma(1)}x_{\sigma(3)}) & \\
+b_3x_{\sigma(3)}(x_{\sigma(2)}x_{\sigma(1)})+b_4x_{\sigma(1)}(x_{\sigma(3)}x_{\sigma(2)})+b_5x_{\sigma(2)}(x_{\sigma(3)}x_{\sigma(1)})+b_6x_{\sigma(3)}(x_{\sigma(1)}x_{\sigma(2)})&=0.
\end{array}
$$
This relation can be also writes
$$\sigma \circ v_1((x_1x_2)x_3)+\sigma \circ v_2(x_1(x_2x_3))=0$$
with $$\sigma \circ (r_1Id+r_2\tau_{1,2}+r_3\tau_{1,3}+r_4\tau_{2,3}+r_5c+r_6c^2)=r_1\sigma+r_2\sigma \circ\tau_{1,2}+r_3\sigma \circ\tau_{1,3}+r_4\sigma \circ\tau_{2,3}+r_5\sigma \circ c+r_6\sigma \circ c^2$$
which corresponds to the pair of vectors $(\sigma \circ v_1,\sigma \circ v_2).$
So we can associate to the relation (\ref{mu}) two matrices $E$ and $F$  whose lines are the components of $\sigma v_1$ and $\sigma  v_2$ related to $\mathcal{B}$, that is 
$$E=\begin{pmatrix}
      a_1&a_2 & a_3 & a_4 &a_5 & a_6    \\
      a_2&a_1 & a_6 & a_5 &a_4 & a_3    \\
        a_3&a_5 & a_1 & a_6 &a_2 & a_4    \\
          a_4&a_6 & a_5 & a_1 &a_3 & a_2    \\
            a_6&a_4 & a_2 & a_3 &a_1 & a_5    \\
              a_5&a_3 & a_4 & a_2 &a_6 & a_1    \\
\end{pmatrix}
\ \ F=\begin{pmatrix}
      b_1&b_2 & b_3 & b_4 &b_5 & b_6    \\
      b_2&b_1 & b_6 & b_5 &b_4 & b_3    \\
        b_3&b_5 & b_1 & b_6 &b_2 & b_4    \\
          b_4&b_6 & b_5 & b_1 &b_3 & b_2    \\
            b_6&b_4 & b_2 & b_3 &b_1 & b_5    \\
              b_5&b_3 & b_4 & b_2 &b_6 & b_1    \\

\end{pmatrix}.
$$
Now we can consider  a linear combination of these relations associated to the pairs of vectors
$(\sigma \circ v_1,\sigma \circ v_2)$
that is which corresponds to the vector
$$((u_1Id+u_2\tau_{1,2}+u_3\tau_{1,3}+u_4\tau_{2,3}+u_5c+u_6c^2)\circ v_1,(u_1Id+u_2\tau_{1,2}+u_3\tau_{1,3}+u_4\tau_{2,3}+u_5c+u_6c^2)\circ v_2).$$
All the coefficients of this relation written in the order associated to the basis $\mathcal{B}$, that is the coefficients of $(x_1x_2)x_3, (x_2x_1)x_3, (x_3x_2)x_1, (x_1x_3)x_2, (x_2x_3)x_1, (x_3x_1)x_2, x_1(x_2x_3),
x_2(x_1x_3),$ 

\noindent $x_3(x_2x_1), x_1(x_3x_2), x_2(x_3x_1), x_3(x_1x_2$ are given by the matricial products
$$AU,BU$$ where $U$ is the column matrix corresponding to the vector $(u_1,u_2,u_3,u_4,u_5,u_6)$ and $A$ and $B$ the matrices

$$A=\begin{pmatrix}
      a_1&a_2 & a_3 & a_4 &a_6 & a_5    \\
      a_2&a_1 & a_5 & a_6 &a_4 & a_3    \\
        a_3&a_6 & a_1 & a_5 &a_2 & a_4    \\
          a_4&a_5 & a_6 & a_1 &a_3 & a_2    \\
            a_5&a_4 & a_2 & a_3 &a_1 & a_6    \\
              a_6&a_3 & a_4 & a_2 &a_5 & a_1    \\
\end{pmatrix}
\ \ B=\begin{pmatrix}
      b_1&b_2 & b_3 & b_4 &b_6 & b_5    \\
      b_2&b_1 & b_5 & b_6 &b_4 & b_3    \\
        b_3&b_6 & b_1 & b_5 &b_2 & b_4    \\
          b_4&b_5 & b_6 & b_1 &b_3 & b_2    \\
            b_5&b_4 & b_2 & b_3 &b_1 & b_6    \\
              b_6&b_3 & b_4 & b_2 &b_5 & b_1    \\
\end{pmatrix}.
$$
Let us note that $A=^tE$, $B=^tF.$
\subsection{Nonassociativity and depolarization}
Assume that $\mu$ is the depolarization of the pair $(\bu,[-,-])$ of a commutative multiplication and a anticommutative multiplication without other hypothesis on these multiplication $\bu$ and $[-,-]$. It is clear that this simple hypothesis on $\mu$ does not generate any relation of the type (\ref{mu}). 

Conversely, consider a nonassociative multiplication satisfying the relation (\ref{mu}). It is the polarization of the couple $(\bu,[-,-])$ with $x\bu y=xy+yx$ and $[x,y]=xy-yx$ (recall that $xy =\mu(x,y)$). What is the quadratic relation between $\bu$ and ${-,-}$ generates by (\ref{mu}) and which restores this identity by depolarization? We have
$$\begin{array}{ll}
     (xy)z =& (x\bu y)\bu z+[x,y]\bu z + [x\bu y,z]+[[x,y],z]   \\
      x(yz)=& x \bu (y \bu z)+x \bu [y,z]=[x,y\bu z]+[x, [y,z]]  
\end{array}
$$
We deduce the relation
\begin{equation}
\label{demu}
\begin{array}{l}
\lambda_1(x_1 \bu x_2)\bu x_3 +\lambda_2(x_2 \bu x_3)\bu x_1 +
\lambda_3(x_1 \bu x_3)\bu x_2 +\lambda_4 [x_1 \bu x_2, x_3 ]+\lambda_5 [x_2 \bu x_3, x_1 ]\\
+\lambda_6 [x_1 \bu x_3, x_2 ]
+\lambda_7[x_1, x_2]\bu x_3 +\lambda_8[x_3 , x_2]\bu x_1 +
\lambda_9[x_1 , x_3]\bu x_2 
+\lambda_{10}[[x_1, x_2], x_3]\\ +\lambda_{11}[[x_3 , x_2], x_1] +
\lambda_{12}[[x_1, x_3], x_2]=0
\end{array}
\end{equation} 
with
$$
\begin{array}{ll}
     \lambda_1= &  a_1+a_2+b_3+b_6 \\
      \lambda_2= &  a_3+a_5+b_1+b_4  \\
       \lambda_3	= &  a_4+a_6+b_2+b_5  \\ 
        \lambda_4= &  a_1+a_2-b_3-b_6    \\
         \lambda_5= &  a_3+a_5-b_1-b_4  \\
          \lambda_6= &   a_4+a_6-b_2-b_5  \\
\end{array}
\ \ \ 
\begin{array}{ll}
     \lambda_7= &   a_1-a_2-b_3+b_6   \\
      \lambda_8= &   a_3-a_5-b_1+b_4 \\
       \lambda_9	= &   a_4-a_6-b_2+b_5 \\ 
        \lambda_{10}= &  a_1-a_2+b_3-b_6    \\
         \lambda_{11}= & a_3-a_5+b_1-b_4   \\
          \lambda_{12}= &  a_4-a_6+b_2-b_5   \\
\end{array}
$$
Since the passage from the variables $a_i$ and $b_j$ to the variables $ \lambda_i$ is reversible, and like any quadratic relations linking the multiplications $\bu$ and $[-,-]$
can be written
$$
\begin{array}{l}
v_1 ((x_1 \bu x_2)\bu x_3 + [x_1 \bu x_2, x_3 ]+[x_1, x_2]\bu x_3 +[[x_1, x_2], x_3])+\\
v_2(( x_1 \bu (x_2\bu x_3) + [x_1, x_2\bu x_3 ]+x_1\bu [ x_2] x_3 ]+[x_1,[ x_2] x_3]])=0,
\end{array}
$$
with
$$v_1=a_1Id+a_2\tau_{12}+a_3\tau_{13}+a_4\tau_{23}+a_5c+a_6c^2 \ \ {\rm and} \ \ v_2=b_1Id+b_2\tau_{12}+b_3\tau_{13}+b_4\tau_{23}+b_5c+b_6c^2$$
we deduce
\begin{proposition}
Any nonassociative multiplication on $V$ satisfying a nonassociativity relation  (\ref{mu}) 
is the depolarization of the couple $(\bu,[-,-])$ which are linked by the relation (\ref{demu}).
\end{proposition}

\subsection{Depolarization and distributive laws}

In this section, we consider two multiplications on $V$, $\bu$ which is commutative and $[-,-]$ which is anticommutative and $\mu$ the depolarization of the pair $(\bu,[-,-] )$. We will say more generally that $\mu$ is the depolarization of a pair (commutative, anticommutative).

A distributive law between these two multiplications  $(\bu,[-,-] )$. is given by a relation of the type
\begin{equation}\label{dis1}\alpha_1 x_1 \bu [x_2,x_3]+\alpha_2 x_2 \bu [x_3,x_1] +\alpha_3 x_3 \bu [x_1,x_2]+\beta_1[x_1 \bu x_2,x_3]+\beta_2[x_2\bu x_3,x_1]+\beta_3 [x_3 \bu x_1,x_2]=0.
\end{equation}
This identity is therefore a special case of (\ref{demu})  corresponding to
$$\lambda_1=\lambda_2=\lambda_3=\lambda_{10}=\lambda_{11}=\lambda_{12}=0$$
and
$$\lambda_4=\beta_1,\lambda_5=\beta_2,\lambda_6=\beta_3,\lambda_7=\alpha_3,\lambda_8=-\alpha_1,\lambda_9=-\alpha_2.$$
We deduce that the depolarization $\mu$ of the couple $(\bu,[-,-])$ satisfies the relation
\begin{equation}
\label{dis2}
 \begin{array}{ll}\rho_1(x_1x_2)x_3+\rho_2(x_2x_1)x_3+\rho_3(x_3x_2)x_1+\rho_4(x_1x_3)x_2+\rho_5(x_2x_3)x_1+\rho_6(x_3x_1)x_2 & \\
-\rho_3x_1(x_2x_3)-\rho_6x_2(x_1x_3)-\rho_1x_3(x_2x_1)-\rho_5x_1(x_3x_2)-\rho_4x_2(x_3x_1)-\rho_2x_3(x_1x_2)&=0
\end{array}        
\end{equation}
with
$$
\begin{array}{l}
     \rho_1=  \alpha_3+\beta_1  \\
      \rho_2=  -\alpha_3+\beta_1  \\
       \rho_3=  -\alpha_1+\beta_2  \\
        \rho_4= - \alpha_2+\beta_3  \\
         \rho_5=  \alpha_1+\beta_2  \\
          \rho_6=  \alpha_2+\beta_3  \\
          \end{array}
         \ \ \ {\rm or \ \ }
\begin{array}{l}
   \alpha_1=  \frac{1}{2}(\rho_5-\rho_3)  \\
      \alpha_2= \frac{1}{2}(\rho_6-\rho_4)  \\
       \alpha_3= \frac{1}{2}(\rho_1-\rho_2)  \\
       \beta_1= \frac{1}{2}(\rho_1+\rho_2)  \\
        \beta_2= \frac{1}{2}(\rho_3+\rho_5)   \\
        \beta_3=  \frac{1}{2}(\rho_4+\rho_6)   \\
          \end{array}         
$$         
          
We therefore associate to a quadratic relation defined by a distributive law the vectors
$$W_1=^t(\rho_1,\rho_2,\rho_3,\rho_4,\rho_5,\rho_6), \ \ W_2=^t(-\rho_3,-\rho_6,-\rho_1,-\rho_5,-\rho_4,-\rho_2).$$

Conversely, if we have a multiplication $\mu$ which satisfies a quadratic relation associated with the pair of vectors $W_1$ and $W_2$, it is the depolarization of the pair $(\bu,[-,-])$ and these multiplications are associated with the distributive law (\ref{dis1}).

\noindent{\bf Examples.}
\begin{enumerate}
  \item  {\bf The Leibniz relation}.  It is a distributive law between $\bu$ and $[-,-]$ given by
  $$[x_1 \bu x_2,  x_3] -x_1 \bu [x_2,x_3] -  [x_1, x_3]\bu x_2=0$$
If $\mu(x,y)=xy$ is a depolarization of the pair $(\bu,[-,-])$, we have the equivalent relation:
\begin{equation}
\label{Lei}\begin{array}{l}
 (x_1  x_2)  x_3 + (x_2x_1) x_3 + (x_3x_2)  x_1 - (x_1  x_3)  x_2-(x_2  x_3) x_1+ (x_3 x_1) x_2 - x_1 (x_2  x_3) - x_2  (x_1  x_3) \\
-x_3  (x_2  x_1) + x_1  (x_3  x_2) + x_2 (x_3  x_1) - x_3  (x_1 x_2) = 0.
 \end{array}
\end{equation}
In this case
$$W_1=(1,1,1,-1,-1,1), \ \ {\rm and } \ \ W_2=(-1,-1,-1,1,1,-1).$$

\medskip

\noindent{\bf Remarks.} 1. The problem we will have to solve is to find which nonassociativity relations described in (\ref{mu}) which by linear combinations of this relation and relations deduced by transposition gives the desired distributive law, for example in the case of the identity of Leibniz gives identity (\ref{Lei}). To do this, we need to determine the matrices $A$ and $B$ that have a vector $U$, such as
$$AU=W_1, \ BU=W_2.$$
For example, if we take $B=Id$, then $U=-W_2=W_1$ and $AW_1=W_1$ gives
$$a_3=1-a_1,a_5=-a_2,a_6=-a_4$$ and $\mu$ have to satisfy the axiom
$$
\begin{array}{l}
a_1(x_1x_2)x_3+a_2(x_2x_1)x_3+(1-a_1)(x_3x_2)x_1+a_4(x_1x_3)x_2-a_2(x_2x_3)x_1-a_4(x_3x_1)x_2 \\
-x_1(x_2x_3)=0.
\end{array}$$

2. We can also be interested in knowing which properties are deduced from the relation (\ref{Lei}). This is equivalent to determining the  images $W'_1,W'_2$ of the matrices $A$ and $B$  associated with this  relation:
$$A U=W'_1, \ BU=W'_2$$ with 
$$A=-B=\begin{pmatrix}
   1 & 1 & 1 &-1 & -1 &1   \\
    1	 & 1 & 1 &-1 & -1 &1   \\
      1 & -1 & 1 &1&1 &-1   \\
 -1 & 1 & -1  &1 & 1 &1   \\
  1 & -1 & 1  &1& 1&-1   \\
 -1 & 1 & -1  &1 & 1 &1   \\
\end{pmatrix}
$$
Since $B=-A$, $W'_2=-W'_1$. The rank of $A$ is $3$ and its image is generated by the columns $C_1(A),C_2(A),C_4(A)$ then $W_1'$ is the transposed of
$$(\alpha_1+\alpha_2-\alpha_4,\alpha_1+\alpha_2-\alpha_4,\alpha_1-\alpha_2+\alpha_4,-\alpha_1+\alpha_2+\alpha_4,\alpha_1-\alpha_2+\alpha_4,-\alpha_1+\alpha_2+\alpha_4)$$ that is 
$$W'_2=^t(\rho_1,\rho_1,\rho_2,\rho_3,\rho_2,\rho_3).$$
Consequence. If the multiplication $\mu$ satisfies the relation (\ref{Lei}), then its associator $\mathcal{A}_\mu$ defined by $\mathcal{A}_\mu(x_1,x_2,x_3)=(x_1x_2)x_3-x_1(x_2x_3)$ satisfies the relation
$$\mathcal{A}_\mu(x_1,x_2,x_3)+\mathcal{A}_\mu(x_2,x_1,x_3)=0$$
that is it is an antiPreLie algebra.

  \item {\bf The transposed Leibniz relation}.  This distributive law is given by
  $$2x_3 \bu [x_1, x_2] = [x_3 \bu x_1, x_2] + [x_1, x_3 \bu x_2].$$
If $\mu(x,y)=xy$ is a depolarization of the pair $(\bu,[-,-])$, we have the equivalent relation:
\begin{equation}
\label{TLei}
\begin{array}{l}
      2(x_1x_2)x_3-2(x_2x_1)x_3+(x_3x_2)x_1-(x_1x_3)x_2+(x_2x_3)x_1-(x_3x_1)x_2  \\
     -x_1(x_2x_3)+x_2(x_1x_3)-2x_3(x_2x_1)-x_1(x_3x_2)+x_2(x_3x_1)+2x_3(x_1x_2)=0  
\end{array}
\end{equation}
In this case we have 
$$W_1=(2,-2,1,-1,1,-1), \ \ {\rm and } \ \ W_2=(-1,1,-2,-1,1,2).$$
\medskip

\noindent{\bf Remarks}
1.  As in the previous example, we can be interested in nonassociative relations whose identity above is a consequence. We have seen hat we have to find $A$ and $B$ such as there exists $U$ with $AU=W_1, \ BU=W_2$.  Take for example $B=-Id$. In this case $U=^t(1,-1,2,1,-1,-2)$.  Then $AU=W_1$ is equivalent to
We have
$$\left\{
\begin{array}{l}
a_1=\frac{7}{3}a_4-2a_5+\frac{2}{3}a_6\\
a_2=2a_4-2a_5+a_6\\
a_3=-\frac{2}{3}a_4+a_5+\frac{2}{3}a_6+1
\end{array}
\right.
$$
and we obtain the general relation of nonassociativitivy implying (\ref{TLei}):
$$
\begin{array}{l}
a_4(\frac{7}{3}(x_1x_2)x_3-2(x_2x_1)x_3-\frac{2}{3}(x_3x_2)x_1+(x_1x_3)x_2)   \\
 +a_5(-2(x_1x_2)x_3-2(x_2x_1)x_3+(x_3x_2)x_1+(x_2x_3)x_1)     \\
 +a_6(\frac{2}{3}(x_1x_2)x_3+(x_2x_1)x_3+\frac{2}{3}(x_3x_2)x_1+(x_3x_1)x_2\\
 +(x_3x_2)x_1-x_1(x_2x_3)=0 
\end{array}
$$

\medskip

2. To determine the relations that arise from (\ref{TLei}), just solve the system
$$AU=W'_1, \ BU=W'_2.$$
this means that $W'_1 \in \im (A)$ and $W'_2 \in \im (B).$
Then $W'_1= \alpha_1C_1(A)+\alpha_2C_3(A)+\alpha_4C_4(A), W'_2= \alpha_1C_1(B)+\alpha_2C_3(B)+\alpha_4C_4(B)$ where $C_i(M)$ is the column $i$ of the matrix $M$.  We thus obtain the following relation
$$
\begin{array}{l}
\alpha_1( 2(x_1x_2)x_3-2(x_2x_1)x_3+(x_3x_2)x_1-(x_1x_3)x_2+(x_2x_3)x_1-(x_3x_1)x_2  \\
     -x_1(x_2x_3)+x_2(x_1x_3)-2x_3(x_2x_1)-x_1(x_3x_2)+x_2(x_3x_1)+2x_3(x_1x_2))\\
     +\alpha_2((x_1x_2)x_3-+(x_2x_1)x_3+2(x_3x_2)x_1-(x_1x_3)x_2-2(x_2x_3)x_1-(x_3x_1)x_2  \\
     -2x_1(x_2x_3)+x_2(x_1x_3)-x_3(x_2x_1)+2x_1(x_3x_2)+x_2(x_3x_1)-x_3(x_1x_2))\\
	+\alpha_3(-(x_1x_2)x_3-(x_2x_1)x_3+(x_3x_2)x_1+2(x_1x_3)x_2+(x_2x_3)x_1-2(x_3x_1)x_2  \\
     -x_1(x_2x_3)+2x_2(x_1x_3)+1x_3(x_2x_1)-x_1(x_3x_2)-2x_2(x_3x_1)+x_3(x_1x_2))=0.
 \end{array}
 $$    
\end{enumerate}

\subsection{General case}

Let $\mu$ be a depolarization of a given couple $(\bullet,[,])$ with $\bullet$ a commutative and associative multiplication and $[-,-]$ an anticommutative multiplication (on $V$).  It is assumed that these two elements of this couple are connected by a distributive law
(\ref{dis1}). This is equivalent to say that $\mu$ satisfies the relation (\ref{dis2}). As in the previous examples, we want to determine the nonassociative relations which  must verify $\mu$ and which involve, by linear combinations associated with the action of $\Sigma_3$, the relation (\ref{dis2}).

This problem consists to find the matrices  $A$ and $B$  and a vector $U$ such as 
$$AU=W_1, \ \ BU=W_2$$
with 
$$W_1=(\rho_1,\rho_2,\rho_3,\rho_4,\rho_5,\rho_6), \ \ W_2=(-\rho_3,-\rho_6,-\rho_1,-\rho_5,-\rho_4,-\rho_2)$$
and $\rho_i \in \K.$
Since $^tW_1 \in \im (A)$, then the vectors $$^t(\rho_2,\rho_1,\rho_6,\rho_5,\rho_4,\rho_3),\cdots,^t(\rho_5,\rho_3,\rho_4,\rho_2,\rho_6,\rho_1)$$ are also in $\im (A)$.  If we denote by $r(\rho)$, the rank of these $6$ vectors, we deduce that the rank of $A$ (and also for the rank of $B$) is greater than $r(\rho)$. For example, if $W_1$ and $W_2$ correspond to the Leibniz relation, then $r(\rho)=3$. Then we have to consider that the rank of the matrices $A$ and $B$ are greater than $3$. If these ranks are $3$, then $A$ and $B$ are defined from $W_1$ and $W_2$ and we find again the Leibniz relation. In the examples given in the previous section, we have taken $B$ whose rank is $6$.

\medskip

In \cite{GRNonass}, we have classified these vectors according to the rank   that is the rank of the $\Sigma_3$-module they generate. We recall this list.
\begin{equation}
\label{liste}
\begin{array}{l}
 Rank =1.\left\{
\begin{array}{l}
V_1^1=Id-\tau_{12}-\tau_{13}-\tau_{23}+c+c^2, \\
V_1^2=Id+\tau_{12}+\tau_{13}+\tau_{23}+c+c^2,
\end{array}
\right.\\
\medskip
Rank=2
\left\{
\begin{array}{l}
V_2^1 = (b_1, -b_1, b_1 + b_5, -b_5, b_5, -b_1-b_5).\\
  V_2^2=(b_1,b_2,b_2,b_2,b_1,b_1), \ \ b_2 \neq \pm b_1.
\end{array}
\right.\\
\medskip
Rank=3
\left\{
\begin{array}{l}
V_3^1=(1 , t,0,-1,0,-t), \ \ t \neq 1,\\
V_3^2=(Id,-1,0,-2, 2,0),\\
V_3^3=(-2,0,-(2 + t), t -1,-(1 + t), t)
\end{array}
\right.\\
\medskip
Rank=4
\left\{
\begin{array}{l}
V_4^1=(2, 1 + t, 1,0,1,1-t), \ \ \ t \neq 1,\\
V_4^2= (2,1 ,0,1,1,1),\\
V_4^3=(2,0,1,-1 ,3,1)\\
V_4^4=(1,0,\alpha,-\alpha,\beta,-1-\beta), \ \alpha^2 \neq 1+\beta+\beta^2
\end{array}
\right.\\
\medskip
Rank=5
\left\{
\begin{array}{l}
V_5^1=(2, -1,-1, -1,1,0),\\
V_5^2=(2,1,1,1,1,0)
\end{array}
\right.\\
\medskip
Rank=6
\left\{
\begin{array}{l}
V_6^1=(1,0,0,0,0,0).
\end{array}
\right.\\
\end{array}
\end{equation}

\noindent {\bf Remark.} In \cite{GRNonass}, the vector $V_4^4$ which corresponds to a $\Sigma_3$-module of dimension $4$ which is the sum of two irreducible modules of dimension $2$ is not very well defined.

For each of the cases listed above, the $B$ matrix is well defined.   We can choose as matrix $B$ the matrix associated with one of these vectors, the only condition is on the rank of $B$ (or $V$) which must be greater or equal to that of $\rho$ and that $W_2 \in \im (B)$. 

\noindent{\bf Example.} Assume that $r(\rho)=1$. In this case $W_2=(\rho_1,-\rho_1,-\rho_1,-\rho_1,\rho_1,\rho_1)$. If we choose $B$ such as $rank(B)=1$, then $B$ is generated by $V^1_1$. Since the choice of $B$ is not determinative, we can always assume that 
$rank(A) \leq rank(B).$  But $W_1=W_2$ and $A$ is also generated by $V_1^1$. In this case $\mu$ satisfies the axiom
$$\begin{array}{l}(x_1x_2)x_3-(x_2x_1)x_3-(x_3x_2)x_1-(x_1x_3)x_2+(x_2x_3)x_1+(x_3x_1)x_2  \\
+x_1(x_2x_3)-x_2(x_1x_3)-x_3(x_2x_1)-x_1(x_3x_2)+x_2(x_3x_1)+x_3(x_1x_2)=0.
\end{array} 
$$ and
the distributivity law 
\begin{equation}
\label{R11}
x_1 \bu [x_2,x_3]+x_2\bu [x_3,x_1]+x_3 \bu [x_1,x_2]=0.
\end{equation}
Let’s now choose for $B$ a matrix of rank $6$. From the list above, $B=Id$. Then $U=W_2$ and, since $W_1=W_2$ the matrix $A$ satisfies $AW_2=W_2$. This is equivalent to have $a_1-a_2-a_3-a_4+a_5+a_6=1$. In this (maximal) case, the relation of nonassociativity of $\mu$ is
$$
\begin{array}{ll}(x_1x_2)x_3+a_2((x_2x_1)x_3+(x_1x_2)x_3)+a_3((x_3x_2)x_1+(x_1x_2)x_3)+a_4((x_1x_3)x_2+(x_1x_2)x_3)\\
+a_5((x_2x_3)x_1-(x_1x_2)x_3+a_6((x_3x_1)x_2 -(x_1x_2)x_3)+x_1(x_2x_3)=0.
\end{array}
$$

\medskip

In the general case, a procedure to determine the nonassociative laws that are the depolarizations of a couple (commutative associative, anticommutative), would consist in starting by determining the "maximal" relations that is to say associated with matrices $B$ of maximum rank. We can then look at the $B$ matrices associated with the $V$ vectors of smaller rank whose orbit is not contained in those of the $V$ vectors already chosen.

For example, take $rank(B)=6$ that is $B=Id$.  Then $BU=W_2$ gives $U=W_2$ and $AU=W_1$ is reduced to $AW_2=W_1$. Since
$$W_1=(\rho_1,\rho_2,\rho_3,\rho_4,\rho_5,\rho_6), \ \ W_2=(-\rho_3,-\rho_6,-\rho_1,-\rho_5,-\rho_4,-\rho_2)$$
then $AW_2=W_1$ is equivalent to 
$$a_1\rho_3+a_2\rho_6+(a_3+1)\rho_1+a_4\rho_5+a_5\rho_2+a_6\rho_4=0$$
We can always assume that $\rho_1 \neq 0$ more precisely $\rho_1=1$. This gives
$$a_3=-a_1\rho_3-a_2\rho_6-a_4\rho_5-a_5\rho_2-a_6\rho_4-1$$
and the multiplication $\mu$ satisfies the following axiom
\begin{equation}
\label{disgene}
\begin{array}{ll}((x_1x_2)x_3-\rho_1(x_3x_2)x_1)+a_2((x_2x_1)x_3-\rho_6(x_3x_2)x_1)-(x_3x_2)x_1\\
+a_4((x_1x_3)x_2-\rho_5(x_3x_2)x_1)
+a_5((x_2x_3)x_1-\rho_2(x_3x_2)x_1)+a_6((x_3x_1)x_2 -\rho_4(x_3x_2)x_1)\\+x_1(x_2x_3)=0.
\end{array}
\end{equation}


\section{Depolarization of the pair (Lie,associative commutative)}
In this section, we consider two multiplications on $V$:
\begin{enumerate}
  \item $[x,y]$ which is a Lie multiplication (or Lie bracket),
  \item $x \bu y$ which is associative and commutative.
\end{enumerate}
The purpose of this section is to describe all multiplications $\mu(x,y)=xy$ on $V$ that are depolarizations of the pair $([-,-],\bu)$. 
As we said in the introduction, we assume that $ \mu$ satisfies a quadratic relation (\ref{mu}):
$$
\begin{array}{ll}a_1(x_1x_2)x_3+a_2(x_2x_1)x_3+a_3(x_3x_2)x_1+a_4(x_1x_3)x_2+a_5(x_2x_3)x_1+a_6(x_3x_1)x_2 & \\
+b_1x_1(x_2x_3)+b_2x_2(x_1x_3)+b_3x_3(x_2x_1)+b_4x_1(x_3x_2)+b_5x_2(x_3x_1)+b_6x_3(x_1x_2)&=0.
\end{array}
$$
\subsection{The Jacobi condition}
The Jacobi identity of $[-;-]$ gives
\begin{equation}\label{jac}
\begin{array}{ll}(x_1x_2)x_3-(x_2x_1)x_3-(x_3x_2)x_1-(x_1x_3)x_2+(x_2x_3)x_1+(x_3x_1)x_2 & \\
-x_1(x_2x_3)+x_2(x_1x_3)+x_3(x_2x_1)+x_1(x_3x_2)-x_2(x_3x_1)-x_3(x_1x_2)&=0.
\end{array}
\end{equation}
that is $\mu$ is Lie admissible. Let $A$ and $B$ the matrices associated with $\mu$. Then $\mu$ is Lie admissible iff 
$\lambda_1=a_1-a_2-a_3-a_4+a_5+a_6 $ and $\lambda_2=b_1-b_2-b_3-b_4+b_5+b_6$ which are eigenvalues of $A$ and $B$ satisfy
$$\lambda_1=-\lambda_1 \neq 0.$$
In fact this is the  necessary and sufficient conditions to have a solution to the linear system $AU=^t(1,-1,-1,-1,1,1), BU=^t(-1,1,1,1,-1,-1)$ (see also \cite{GRNonass}).
So, $\mu$ satisfy the relation
\begin{equation}\label{muJac}
\begin{array}{ll}a_1(x_1x_2)x_3+a_2(x_2x_1)x_3+a_3(x_3x_2)x_1+a_4(x_1x_3)x_2+a_5(x_2x_3)x_1\\
+(\lambda-a_1+a_2+a_3+a_4-a_5)(x_3x_1)x_2 & \\
+b_1x_1(x_2x_3)+b_2x_2(x_1x_3)+b_3x_3(x_2x_1)+b_4x_1(x_3x_2)+b_5x_2(x_3x_1)\\+(-\lambda-b_1+b_2+b_3+b_4-b_5)x_3(x_1x_2)&=0.
\end{array}
\end{equation}

\subsection{The associativity of $\bullet$}
The associativity of the commutative multiplication $\bu$ is equivalent to
\begin{equation}\label{ass}
(x_1x_2)x_3+(x_2x_1)x_3-(x_3x_2)x_1-(x_2x_3)x_1
-x_1(x_2x_3)+x_3(x_2x_1)-x_1(x_3x_2)+x_3(x_1x_2)=0.
\end{equation}
In this case, we are led to find the matrices $A$ and $B$ and a vector $U$ such as
$$AU=W_3=^t(1,1,-1,0,-1,0), \ BU=W_4=^t(-1,0,1,-1,0,1).$$
In particular $W_4 \in \im (B)$ that implies that the vectors $^t(0,-1,1,0,-1,1), \ ^t(-1,1,0,-1,1,0)$ are also in this vector space.  Since these three vectors are dependent, we have $rank(B) \geq 2$ and also for $A$.  The minimal solution, that is the solution corresponding to the minimal rank is given by $B$ (resp.$A$) generated by $W_4$ (resp $W_3$), but in this case we find again (\ref{ass}).  Let us note that the vector $W_4$ generates the irreducible $\Sigma_3$-module of dimension $2$ (see also \cite{GRNonass}, Proposition 5) . We will therefore look for the maximum solution, that is to say the one corresponding to the rank equal to $6$. In this case, we can take $B=Id$ implying $U=^t(-1,0,1,-1,0,1)$ and $AU=^t(1,1,-1,0,-1,0)$ is equivalent to 
$$a_5=1+a_1-a_3+a_4, \ a_6=a_1-a_2+a_4.$$
In this case $\mu$ have to satisfy the axiom
\begin{equation}\label{muAss}
\begin{array}{l}
a_1((x_1x_2)x_3+(x_2x_3)x_1+(x_3x_1)x_2)+a_2((x_2x_1)x_3-(x_3x_1)x_2))+a_3((x_3x_2)x_1-(x_2x_3)x_1)\\
+a_4((x_1x_3)x_2+(x_2x_3)x_1+(x_3x_1)x_2)+(x_2x_3)x_1+x_1(x_2x_3)=0
\end{array}
\end{equation}

\medskip

\noindent{\bf Remarks.} We have seen that $rank(B) \geq 2$. If this rank is $2$, we find again the relation giving the associativity of $\bu$. We have also seen that the vector $W_4=^t(-1,0,1,-1,0,1)$ generates the irreducible $\Sigma-3$-module of dimension $2$. Then $\mu$ is not Lie admissible nor $3$-power associative because these two properties are associated with the two irreducible $\Sigma_3$-modules of dimension $1$. This implies also that the matrix $B$ which defines $\mu$ cannot be of rank $3$ or $5$ (\ref{liste}). When the rank is $4$,  we have to consider the vectors $V_4^1(2, 1 + t, 1,0,1,1-t), \ \ \ t \neq 1,
V_4^2= (2,1 ,0,1,1,1),
V_4^3=(2,0,1,-1 ,3,1), V_4^4=(1,0,\alpha,-\alpha,\beta,-1-\beta), \ \alpha^2 \neq 1+\beta+\beta^2$,  The modules corresponding to the first three vectors contain an irreducible module of dimension $1$ and cannot be suitable. Then we only need to consider the matrices associated with vector $V_4^4$. In this case the linear system $BU=W_4$ with $B$ defined by $V_4^4$ gives 
$U=^t(\frac{-1-\alpha}{\Delta},\frac{-\beta}{\Delta},\frac{1+\alpha}{\Delta},0,\frac{\beta}{\Delta},0)$ with $\Delta =1=\beta+\beta^2-\alpha^2$ which is, by hypothesis, not equal to $0$. Then $AU=W_3$ gives $(a_1,a_2,a_3,a_4,a_5,a_6)=$
$$
\begin{array}{l}
 (a_1,a_2,a_3,a_4,a_5,a_6)= ^t(\alpha,-1-\beta,1,\beta,0,-\alpha) \\
\medskip
 \ds +s_1.^t(\frac{-1 -\alpha}{\beta}, 1, \frac{-1 - \alpha}{\beta}, 0, 0, 1)+s_2.^t(\frac{1 +\alpha}{1 + \alpha + \beta},\frac{1 +\alpha}{1 + \alpha + \beta}, 1, 0, 1, 0)\\
\medskip
\ds+s_3.\ds ^t(\frac{1+2a+a^2+b+ab+b^2}{b(1+a+b)},-\frac{1 +\alpha}{1 + \alpha + \beta},\frac{1+\alpha}{\beta},1,0,0).
\end{array}
$$

\subsection{The pair ($\bu$, $[-,-]$)}
Now suppose that $\mu$ is the depolarization of the pair $(\bullet, [-,-])$ where $\bullet$ is  a commutative associative multiplication
and $[-,-]$ a Lie bracket.  We have highlighted two axioms that $\mu$ must verify in order for $ bu$ to be associative. As the second, the one associated with the vector $V_4^4$ is incompatible with the Lie admissibility, we will assume  that $\mu$ satisfies (\ref{muAss}). In this case, the Lie admissibility condition is equivalent to
$$ a_4=-2-3a_1+2a_2+2a_3, \ a_5=-1-2a_1+2a_2+a_3, \ a_6=-2-2a_1+a_2+2a_3$$
and we have
\begin{proposition}
Let $\mu$ be the depolarization of a couple (associative commutative, Lie bracket).  Then $\mu$ have to satisfy the axiom
\begin{equation}
\label{JacAss}
\begin{array}{l}
a_1(x_1x_2)x_3+a_2(x_2x_1)x_3+a_3(x_3x_2)x_1+(-2-3a_1+2a_2+2a_3)(x_1x_3)x_2\\
+(-1-2a_1+2a_2+a_3)(x_2x_3)x_1+(-2-2a_1+a_2+2a_3)(x_3x_1)x_2+x_1(x_2x_3)=0
\end{array}
\end{equation}
or any $\Sigma_3$-linear combination.
\end{proposition}
If $\Phi(x_1,x_2,x_3)$ is a linear combination of the elements $(x_Ix_j)x_k$ and $x_i(x_jx_k)$, by a $\Sigma_3$-linear combination of $\Phi$ we means $\sum u_\sigma\Phi \circ \sigma$.


\section{Depolarization of a triple (Lie,AssComm,distrib)}
Let us return to the general case before concentrating on the case of Poisson and  transposed Poisson algebras. The aim is to depolarize a triplet given by a Lie bracket, an associative and commutative multiplication and a law of distributivity between them. We therefore seek to determine the coefficients of the relations (\ref{mu}) so that by depolarization we find our triplet.  We have previously determined the general form of (\ref{mu}) so that the multiplications are Lie and associative commutative. Let us now seek the conditions for a distributive law to operate between these two multiplications.

A general distributive law has been written on the following form (\ref{dis2}):
$$
 \begin{array}{ll}\rho_1(x_1x_2)x_3+\rho_2(x_2x_1)x_3+\rho_3(x_3x_2)x_1+\rho_4(x_1x_3)x_2+\rho_5(x_2x_3)x_1+\rho_6(x_3x_1)x_2 & \\
-\rho_3x_1(x_2x_3)-\rho_6x_2(x_1x_3)-\rho_1x_3(x_2x_1)-\rho_5x_1(x_3x_2)-\rho_4x_2(x_3x_1)-\rho_2x_3(x_1x_2)&=0
\end{array}   
$$
and our problem consists to solve the linear systems
$$AU=^t(\rho_1,\rho_2,\rho_3,\rho_4,\rho_5,\rho_6),\ \  BU=^t(-\rho_3,-\rho_6,-\rho_1,-\rho_5,-\rho_4,-\rho_2)$$
where $A$ and $B$ are the matrices associated with the axiom (\ref{JacAss}). Since $B=Id$, then
$^tU=^t(-\rho_3,-\rho_6,-\rho_1,-\rho_5,-\rho_4,-\rho_2).$ We have to solve the linear system
$$A.\begin{pmatrix}
      -\rho_3\\
      -\rho_6 \\
      -\rho_1 \\
      -\rho_5 \\
      -\rho_4 \\
      -\rho_2 
\end{pmatrix}
=\begin{pmatrix}
      \rho_1\\
      \rho_2 \\
      \rho_3 \\
      \rho_4 \\
      \rho_5 \\
      \rho_6 
\end{pmatrix}
$$ where $A$ is associated with the vector $(a_1,a_2 , a_3 , -2-3a_1+2a_2+2a_3 ,-2-2a_1+a_2+2a_3, -1-2a_1+2a_2+a_3).$ 
Without any conditions on $A$, we have the general solution
$$\rho_1=-s_1-2s_2-2s_3,\rho_2=s_1,\rho_3=2s_1+2s_2+3s_3,\rho_4=s_2,\rho_5=s_3,,\rho_6=-2s_1-s_2-2s_3$$
that is
$$\rho_1=-\rho_2-2\rho_4-2\rho_5, \ \rho_3=2\rho_2+2\rho_4+3\rho_5, \ \rho_6=-2\rho_2-\rho_4-2\rho_5.$$This means that $\mu$ satisfies also a distributive law
\begin{equation}
\label{disgen}
 \begin{array}{l}
 \rho_2(-(x_1x_2)x_3+(x_2x_1)x_3+2(x_3x_2)x_1-2(x_3x_1)x_2-2x_1(x_2x_3)+2x_2(x_1x_3)+x_3(x_2x_1)\\-x_3(x_1x_2))+\rho_4(-2(x_1x_2)x_3+2(x_3x_2)x_1+(x_1x_3)x_2-(x_3x_1)x_2-2x_1(x_2x_3)+x_2(x_1x_3)\\+2x_3(x_2x_1)-)x_2(x_3x_1))
 +\rho_5(-2(x_1x_2)x_3+3(x_3x_2)x_1+(x_2x_3)x_1-2(x_3x_1)x_2)\\-3x_1(x_2x_3)+2x_2(x_1x_3)+2x_3(x_2x_1)-x_1(x_3x_2)=0
\end{array}        
\end{equation}
that corresponds to
\begin{equation}\label{disgen1}\alpha( x_1 \bu [x_2,x_3]+ x_2 \bu [x_3,x_1] + x_3 \bu [x_1,x_2])+\beta_1[x_1 \bu x_2,x_3]-(\beta_1+\beta_3)[x_2\bu x_3,x_1]+\beta_3 [x_3 \bu x_1,x_2]=0.
\end{equation}

\medskip

The distributivity relation that we have just described is in a way a "universal" relation that is true for any nonassociative law resulting from a depolarization of an associative commutative multiplication and a Lie bracket.  We will now describe distributivity laws based on the $A$ matrix. The linear system $A^t(-\rho_3,-\rho_6, -\rho_1, -\rho_5, -\rho_4 ,-\rho_2 )=^t(\rho_1,\rho_2,\rho_3,\rho_4,\rho_5,\rho_6)$ can be solved with two different views:
\begin{enumerate}
  \item The matrix $A$ is given, we then determine the coefficients $\rho_i$. Note that $A$ depends only on three parameters $a_1,a_2,a_3$.
  \item Coefficients $\rho_i$ are given (as for example for Leibniz relations) and we have to compute the coefficients $a_i$.  This corresponds to solving the system
  $$
  \begin{pmatrix}
     \rho_3 &  \rho_6 &  \rho_1 & \rho_5 &  \rho_2 & \rho_ 4 \\
      \rho_6 &  \rho_3 &  \rho_2 & \rho_4 &  \rho_1 & \rho_ 5 \\ 
       \rho_1 &  \rho_4 &  \rho_3 & \rho_2 &  \rho_5 & \rho_ 6 \\
        \rho_5 &  \rho_2 &  \rho_4 & \rho_3 &  \rho_6 & \rho_ 1 \\
         \rho_4 &  \rho_1 &  \rho_5 & \rho_6 &  \rho_3 & \rho_ 2 \\
          \rho_2 &  \rho_5 &  \rho_6 & \rho_1 &  \rho_4 & \rho_ 3 \\
\end{pmatrix}
\begin{pmatrix}
      a_1    \\
      a_2    \\
      a_3+1    \\
      a_4    \\
       a_5    \\
       a_6    \\
\end{pmatrix}=0
$$
For example, if the rank of the matrix $(\rho_i)$ described in this linear system is of maximal rank $6$, then $a_i=0, i \neq 3$ and $a_3=-1$. But $a_1=a_2=0,a_3=-1$ implies $a_4=-2,a_5=-1,a_6=-4$. Then the distributive law associated with this case cannot be deduced from (\ref{JacAss}).
\end{enumerate}


\section{Depolarization of a Poisson algebra}

The depolarization of Poisson algebras has already been studied in \cite{MRPoisson}. It is used to describe some properties of Poisson algebras \cite{GR,BenaChopp}. We resume here this study in the framework presented in introduction namely, a depolarization of a couple (Lie, commutative associative) and an interpretation of a distributive law represented here by the identity of Leibniz. This section is only to illustrate this depolarization method.

\medskip

Poisson algebras are usually defined as structures with two operations, a commutative
associative one and an anti-commutative one satisfying the Jacobi identity. These
operations are tied up by a distributive law
$$[x \bu y,  z] -x \bu [y, z] -  [x, z]\bu y=0$$
also called the Leibniz condition. If $\mu$ is a depolarization of a Poisson algebra, it is of course a depolarization of the pair (Lie, associative commutative). Thus the Leibniz condition is equivalent to
$$ (x  y)  z + (yx) z + (zy)  x - (x  z)  y-(y  z) x+ (z x) y - x (y  z) - y  (x  z)-z  (y  x) + x  (z  y) + y (z  x) - z  (x y) = 0.$$
With the previous notations, we have 
$$(\rho_1,\rho_2,\rho_3,\rho_4,\rho_5,\rho_6)=(1,1,1,-1,1,-1).$$
Let us consider $\mu$ a depolarization of $\bu$ and $[-,-]$ which satisfies (\ref{JacAss}). Then the Leibniz condition is deduced of (\ref{JacAss}) if and only if the linear system
$$A.^t(-1,-1,-1,1,1,-1)=^t(1,1,1,-1,-1,1)$$
admits a solution. This system is written
$$
\left\{
\begin{array}{l}
   -4a_1+2a_3=4 \\
   2a_1-4a_3=-2 \\
   8a_1-6a_2-4a_3=-6 
\end{array}
\right.
$$
and the solution is $a_1=-1,a_2=-\frac{1}{3},a_3=0$, that implies $a_4=a_5=-a_6=\frac{1}{3}.$

\begin{theorem}\cite{MRPoisson}
The polarization of a Poisson algebra is given by a multiplication $\mu(x,y)=xy$ which satisfies the axiom:
$$3(x y) z +(y x) z -(x z) y - (y z) x + (z x) y-3x(yz)=0.$$
\end{theorem}

\section{Depolarization of transposed Poisson algebras}

\subsection{Definition} \cite{Bai}.  
 Let $V$ be a vector space equipped with two bilinear operations
$$\bu, [- ,- ] :V \otimes V \rightarrow V$$
The triple $(V,\bu, [- ,- ])$ is called a  transposed
Poisson algebra  if $(V,\bu)$ is a commutative associative algebra and $(V, [ -,- ])$ is a
Lie algebra that satisfy the following compatibility condition
\begin{equation}
\label{tp}
2z \bu [x, y] = [z \bu x, y] + [x, z \bu y]
\end{equation}
for any $x,y,z \in V$. The equation (\ref{tp}) is called the transposed Leibniz rule because the roles played by the two
binary operations in the Leibniz rule in a Poisson algebra are switched. 

\medskip

This notion of transposed Poisson algebra has been introduced in \cite{Bai}  and arises naturally from a
Novikov-Poisson algebra by taking the commutator Lie algebra of the Novikov algebra. The most classical example is given taking $V=\mathcal{C}^\infty(\R,\R)$ with the Lie bracket
$$[f,g](x)=f'(x)g(x)-f(x)g'x)$$
and the multiplication $f \bu g =fg.$
In fact 
$$\begin{array}{ll}
    \lbrack fg,h \rbrack + \lbrack  f,gh] \rbrack   & =(f'g+fg')h-h'fg+f'gh-fg'h-fgh'   \\
      &   =2f'gh-2fgh'\\
      &= 2g[f,h]
\end{array}
$$

\subsection{Depolarization of a transposed Poisson algebra}

Since $\bu$ is commutative and associative and $[-,-]$ is a Lie bracket, we have seen that any depolarization of this couple satisfies the axiom (\ref{JacAss}):
$$
\begin{array}{l}
a_1(x_1x_2)x_3+a_2(x_2x_1)x_3+a_3(x_3x_2)x_1+(-2-3a_1+2a_2+2a_3)(x_1x_3)x_2\\
+(-1-2a_1+2a_2+a_3)(x_2x_3)x_1+(-2-2a_1+a_2+2a_3)(x_3x_1)x_2+x_1(x_2x_3)=0.
\end{array}
$$
or any $\Sigma_3$ linear combinations. 
Now, the identity (\ref{tp}) which defines the transposed Leibniz identity is given in (\ref{TLei}):
$$\begin{array}{l}
      2(x_1x_2)x_3-2(x_2x_1)x_3+(x_3x_2)x_1-(x_1x_3)x_2+(x_2x_3)x_1-(x_3x_1)x_2  \\
     -x_1(x_2x_3)+x_2(x_1x_3)-2x_3(x_2x_1)-x_1(x_3x_2)+x_2(x_3x_1)+2x_3(x_1x_2)=0  
\end{array}
$$
Let
$$W_1=(2,-2,1,-1,1,-1), \ \ {\rm and } \ \ W_2=(-1,1,-2,-1,1,2)$$
be the vectors associated with this identity. If $A$ and $B$ are the matrices associated with (\ref{JacAss}), we have to solve the system $AU=W_1,BU=W_2$. Since $B=Id$, then $U=W_2$. The relation $AW_2=W_1$ is equivalent to
$$
\left\{
\begin{array}{l}
 -a_1+a_2=1    \\
 -4a_1+2a_2+2a_3=3\\
 2a_1-2a_3=-3.  
\end{array}
\right.
$$
But this linear system has no solution, that is the relation \ref{JacAss} cannot imply \ref{TLei}).
\begin{proposition}
The depolarization $\mu$ of the couple $(\bu,[-,-]$ which defines a structure of Transposed Poisson algebra on $V$ have to satisfy the  independent axioms
$$
\left\{
 \begin{array}{l}
  2(x_1x_2)x_3-2(x_2x_1)x_3+(x_3x_2)x_1-(x_1x_3)x_2+(x_2x_3)x_1-(x_3x_1)x_2  \\
     -x_1(x_2x_3)+x_2(x_1x_3)-2x_3(x_2x_1)-x_1(x_3x_2)+x_2(x_3x_1)+2x_3(x_1x_2) \\
\\
      (x_1x_2)x_3-(x_2x_1)x_3-(x_3x_2)x_1-(x_1x_3)x_2+(x_2x_3)x_1+(x_3x_1)x_2  \\
     -x_1(x_2x_3)+x_2(x_1x_3)+x_3(x_2x_1)+x_1(x_3x_2)-x_2(x_3x_1)-x_3(x_1x_2) \\
\\
      (x_1x_2)x_3+(x_2x_1)x_3-(x_3x_2)x_1-(x_2x_3)x_1  \\
     -x_1(x_2x_3)+x_3(x_2x_1)-x_1(x_3x_2)+x_3(x_1x_2)  
\end{array}
\right.
$$

\end{proposition}
To say that these  axioms are independent means that we cannot deduce one of these axiom from the second using linear combination of the image by the action of $\Sigma_3$ of this axiom.

\medskip

\noindent{\bf Remark.} Depolarization of a $(a,b,c)$-transposed Poisson algebra.
We can generalize the   transposed Leibniz relation  by considering a relationship of the type
$$a x_3 \bullet [x_1,x_2]-b[x_3 \bullet x_1,x_2]-c[x_1,x_3\bullet x_2]=0$$
with $a \neq 0$ where $[,]$ is an anticommutative multiplication and $\bullet$ a commutative multiplication. This is equivalent to saying that the depolarization $\mu$ of the couple $(\bullet,[-,-])$  satisfies the identity
\begin{equation}
\label{tpabc}
\begin{array}{l}
      a(x_1x_2)x_3-a(x_2x_1)x_3+c(x_3x_2)x_1-b(x_1x_3)x_2+c(x_2x_3)x_1-b(x_3x_1)x_2  \\
     -cx_1(x_2x_3)+bx_2(x_1x_3)-ax_3(x_2x_1)-cx_1(x_3x_2)+bx_2(x_3x_1)+ax_3(x_1x_2)=0  
\end{array}
\end{equation}
or any $\Sigma_3$ linear combinations. 
In this case also, the axiom (\ref{tpabc}) cannot be a consequence of (\ref{JacAss}) which reflect all the depolarization of the couple (Associative commutative, Lie bracket).

\medskip

Depolarization relations allow to find properties of transposed Poisson algebras. Let $A_1,B_1$ and $A_2,B_2$ be the matrices associated respectively to (\ref{JacAss}) and (\ref{TLei}). let us consider the vector  $U=^t(1,-1,-1,-1,1,1)$. We have $A_1U=-U,B_1U=U$ (which ranslate the Lie admissibility of $\mu$) and $A_2=B_2U=4U.$ We deduce that $\mu$ satisfies 
$$
\begin{array}{l}
      (x_1x_2)x_3-(x_2x_1)x_3-(x_3x_2)x_1-(x_1x_3)x_2+(x_2x_3)x_1+(x_3x_1)x_2  \\
     +x_1(x_2x_3)-x_2(x_1x_3)-x_3(x_2x_1)-x_1(x_3x_2)+x_2(x_3x_1)+x_3(x_1x_2)=0  
\end{array}
$$
that is
$$\sum_{\sigma \in \Sigma_3}\mathcal{AA}(x_{\sigma(1)},x_{\sigma(2)},x_{\sigma(3)})=0$$
where $\mathcal{AA}(x_1,x_2,x_3)=(x_1x_2)x_3+x_1(x_2x_3)$ is the antiassociator of $\mu$ (\cite{REAntiass}).
This is equivalent to say
\begin{proposition} (see also \cite{Bai}, theorem 5.)
Any transposed Poisson algebra satisfies also the distributive law
\begin{equation}
\label{AA-A}
x_1 \bullet [x_2,x_3]+x_2 \bullet [x_3,x_1]+x_3 \bullet [x_1,x_2]=0.
\end{equation}
\end{proposition}

\medskip

\noindent{\bf Remarks.} 1. Let $\mu$ be a multiplication which satisfies the identity (\ref{TLei}). Then it satisfies also
$$\sum_{\sigma \in \Sigma_3}\mathcal{AA}(x_{\sigma(1)},x_{\sigma(2)},x_{\sigma(3)})=0.$$
If $\mu$ is Lie admissible, and in this case the second axiom of the depolarization presented in Proposition 4 is satisfied, we have also
$$\sum_{\sigma \in \Sigma_3}\epsilon(\sigma)\mathcal{A}(x_{\sigma(1)},x_{\sigma(2)},x_{\sigma(3)})=0.$$
These two identities are equivalent to
$$
\begin{array}{l}
      (x_1x_2)x_3-(x_2x_1)x_3-(x_3x_2)x_1-(x_1x_3)x_2+(x_2x_3)x_1+(x_3x_1)x_2 =0 \\
     x_1(x_2x_3)-x_2(x_1x_3)-x_3(x_2x_1)-x_1(x_3x_2)+x_2(x_3x_1)+x_3(x_1x_2)=0  
\end{array}
$$
\noindent{\bf Example.} If $V=\mathcal{C}^1(\R,\R)$ the vector space of real functions of class $\mathcal{C}^1$ and
$$f\circ g=f\cdot g, \  \{f,g\}=f'g-fg'$$
then 
$$\mu(f,g)=fg=\frac{1}{2}(f\cdot g +f'g-fg')$$
and this depolarization satisfies
$$(fg)h-(gf)h-(hg)f-(fh)g+(gh)f+(hf)g=0$$
and in this case, the Transposed Leibniz identity implies that the anticommutative multiplication is a Lie bracket. 

2. Let us assume now that $\mu$ which satisfies alwaus (\ref{TLei}) is "associative admissible", that is $\bu$ is associative. In this case the sytem of equations which is construct from (\ref{TLei}) and (\ref{ass}) and with the action of  $\Sigma_3$ is of rank $5$. The matrix in the basis $(x_ix_j)x_k,x_i(x_jx_k)$ in the conventional order is
$$\left(
\begin{array}{ccccccccccccc}
    2  & -2 & 1 & -1 & 1 & -1 & . & -1 & 1 & -2 & -1 & 1 & 2   \\
     1  & 1 & 2 & -1 & -2 & -1 &. & -2 & 1 & -1 &2  & 1 & _1   \\
       -1  & -1 & 1 & 2 & 1 & -2 & .& -1 & 2 &1 & -1 & -2 & 1   \\  
       1  & 1 & 0 & -1 & 0 & -1 & .& 0 & -1 & 1 & 0 & -1 &1   \\
         1  & 1 & -1 & 0 & -1 & 0 & .& -1 & 0 & 1 & -1 & 0 & 1   \\
      \end{array}
      \right)
      $$
      It is easy to see that the multiplication $fg$ described in the previous example satisfies this linear system.

  \subsection{The operad $\mathcal{TP}ois$}

 We recall briefly some facts about operads. A symmetric operad $\mathcal{P}$ is a sequence of vector spaces $\mathcal{P}(n)$, for $n \geq  1$, such that $\mathcal{P}(n)$ is a module over the symmetric group $\Sigma_n$, together with composition maps
$$
\circ_i: \mathcal P(n) \times \mathcal P(m) \to \mathcal P(n+m-1)
$$
satisfying associativity-like conditions:
$$
(f \circ_i g) \circ_j h = f \circ_j (g \circ_{i-j+1} h). 
$$
 An algebra over an operad $\mathcal{P}$ is a vector space $A$ together with maps $\mathcal{P}(n) \otimes A^{\otimes^n} \ra A$ satisfying some  associativity, unitarity and equivariance relations.

The operad for transposed-Poisson algebras
 will be denoted by $\mathcal{TP}ois$.
It is a binary quadratic operad, that is, operad of the form $\p =\Gamma(E)/(R)$, where $\Gamma(E)$ denotes the free operad generated by a
$\Sigma_2$-module $E$ placed in arity 2 and $(R)$ is the operadic ideal generated by a $\Sigma_3$-invariant subspace $R$
of $\Gamma(E)(3).$ To define $\mathcal{TP}ois$, we consider the operadic ideal $R \subset \Gamma(E)(3)$ generated by the vectors 
$$
\left\{
 \begin{array}{l}
  2(x_1x_2)x_3-2(x_2x_1)x_3+(x_3x_2)x_1-(x_1x_3)x_2+(x_2x_3)x_1-(x_3x_1)x_2  \\
     -x_1(x_2x_3)+x_2(x_1x_3)-2x_3(x_2x_1)-x_1(x_3x_2)+x_2(x_3x_1)+2x_3(x_1x_2) \\
\\
      (x_1x_2)x_3-(x_2x_1)x_3-(x_3x_2)x_1-(x_1x_3)x_2+(x_2x_3)x_1+(x_3x_1)x_2  \\
     -x_1(x_2x_3)+x_2(x_1x_3)+x_3(x_2x_1)+x_1(x_3x_2)-x_2(x_3x_1)-x_3(x_1x_2) \\
\\
      (x_1x_2)x_3+(x_2x_1)x_3-(x_3x_2)x_1-(x_2x_3)x_1  \\
     -x_1(x_2x_3)+x_3(x_2x_1)-x_1(x_3x_2)+x_3(x_1x_2)  
\end{array}
\right.
$$

 Then the first terms of $\mathcal{TP}ois=\oplus_{n \geq 1}\mathcal{TP}ois(n)$ are
 $$\mathcal{TP}ois(1)=\K,\mathcal{TP}ois(2)=\Gamma(E)(2)=\K\{x_1x_2,x_2x_1\}.$$ 

 To compute the dimension and the generators of $\mathcal{TP}ois(3)=\Gamma(E)(3)/R$, we consider the three vectors which generate $R$ and their images by the action of $\Sigma_3$. This gives $18$ vectors, but it is easy to see that these vectors are dependent and a basis of this space is given by the vectors whose matrix in the basis $\{(x_{\sigma(1)}x_{\sigma(2)})x_{\sigma(3)},x_{\sigma(1)}(x_{\sigma(2)}x_{\sigma(3)}\}$ with $\sigma=Id,\tau_{12},\tau_{13},\tau_{23},c,c^2$, in this order, is
 $$TP= \left(
 \begin{array}{cccccccccccc}
     2 & -2 & 1 & -1 & 1 & -1 & -1 & 1 & -2 & -1 & 1 & 2   \\
    1 & 1 & 2 & -1 & -2 & -1 & -2 & 1 & -1 & 2 & 1 & -1  \\
    -1 & -1 & 1 & 2 & 1 &-2 & -1 & 2 & 1 & -1 & -2 & 1\\
    1 & -1 & -1 & -1 & 1 & 1 & -1 & 1 & 1 & 1 & -1 & -1\\
    1 & 1 & -1 & 0 & -1 & 0 & -1 & 0 & 1 & -1 & 0 & 1\\
    1 & 1 & 0 & -1 & 0 & -1 & 0 & -1 & 1 & 0 & -1 & 1
\end{array}
\right)
$$
The rank of this matrix is $6$ then
$\dim \mathcal{TP}oisson(3)=6.$

Recall also the definition of the quadratic dual operad $\p^! $ of an operad $\p$. If $\p = \Gamma(E)/(R)$, then there is a scalar product on $\Gamma(E)(3)$ defined by for $\{i,j,k\} = \{i',j',k'\}=\{1,2,3\}$ by
\begin{eqnarray}
\label{pairing}
\left\{
\begin{array}{l}
<(x_i \cdot x_j)\cdot x_k,(x_{i'} \cdot x_{j'})\cdot x_{k'}>=0, \ {\rm if} \ (i,j,k) \neq (i',j',k'), \\
<(x_i \cdot x_j)\cdot x_k,(x_i \cdot x_j)\cdot x_k>={\varepsilon(\sigma)},  
\qquad   {\rm with} \ \sigma =
\left(
\begin{array}{lll}
1 & 2 & 3\\
i &j &k 
\end{array}
\right)
\\
<x_i \cdot (x_j\cdot x_k),x_{i'} \cdot (x_{j'}\cdot x_{k'})>=0, \ {\rm if} \ ( i,j,k) \neq ( i',j',k'), \\
<x_i \cdot (x_j\cdot x_k),x_i \cdot (x_j\cdot x_k)>=-{\varepsilon(\sigma)} 
\qquad {\rm with} \ \sigma =
\left(
\begin{array}{lll}
1 & 2 & 3\\
i &j &k 
\end{array}
\right)
, \\
<(x_i \cdot x_j)\cdot x_k,x_{i'} \cdot (x_{j'}\cdot x_{k'})>=0,
\end{array}
\right.
\end{eqnarray}
where $\varepsilon(\sigma)$ is the signature of $\sigma$.
Let $R^\perp$ be the annihilator of $R$ with respect to this scalar product. Any vector $u=(u_i), i=1,\cdots,12$ which belongs to $R^\perp$ is in the kernel of 
$$DTP= \left(
 \begin{array}{cccccccccccc}
     2 & 2 & -1 & 1 & 1 & -1 &   1 & 1 & -2 & -1 & -1 & -2   \\
    1 & -1 & -2 & 1 & -2 & -1 &    2 & 1 & -1 & 2 & -1 & 1  \\
    -1 & 1 & -1 & -2 & 1 &-2 &   1 & 2 & 1 & -1 & 2 & -1\\
    1 & 1 & 1 & 1 & 1 & 1 &     1 & 1 & 1 & 1 & 1 & 1\\
    1 & -1 & 1 & 0 & -1 & 0 &   1 & 0 & 1 & -1 & 0 & -1\\
    1 & -1 & 0 & 1 & 0 & -1 &   0 & -1 & 1 & 0 & 1 & -1
\end{array}
\right)
$$
Since the rank of this matrix is equal to $6$ and the dual operad $\mathcal{TP}ois^! $ is defined by the operatic ideal $R^\perp$, we have $$\dim \mathcal{TP}ois^! (3)   =6.$$
Since we have also
$$TP.^tDTP=0$$
we can conclude that this operad is self dual that is $\mathcal{TP}ois^! =\mathcal{TP}ois.$ This implies
\begin{proposition}\cite{Bai}
The quadratic operad $\mathcal{TP}ois$ is a Koszul operad.
\end{proposition}

\noindent{\bf Remark.} We could have shown this result by using the theorem of Markl \cite{Ma} which specifies that the operad associated with a distributive law between two Koszul operad of Koszul is also  Koszul.

\medskip

\noindent Recall that this property means that the free $\mathcal{TP}ois$-algebras are Koszul algebras. To illustre this, we can compute the free algebra with one generator $X$. Let us denote by $\mathcal{F}_{TP}(X)=\sum_n\mathcal{F}_{TP}^n(X)$ this algebra. We have
$$\mathcal{F}_{TP}^0(X)=\K, \ \ \mathcal{F}_{TP}^1(X)=\K\{X\}, \ \mathcal{F}_{TP}^2(X)=\K\{X^2\}, \  \mathcal{F}_{TP}^3(X)=\K\{X^2X,XX^2\}$$
For the degree $4$, the quadratic relations which concerns $\mu$ give
$$
\left\{
\begin{array}{l}
4X^2X^2-(x^2X)X-(XX^2)X_X(XX^2)-X(X^2X)=0\\
3(XX^2)X-(X^2X)X-3X(X^2X)+X(XX^2)=0
\end{array}
\right.
$$
and $\dim \mathcal{F}_{TP}^4(X)=3.$


\section{Graduate algebras}
\subsection{Poisson superalgebras}

By a $\K$-super vector space, we mean a $\Z_2$-graded vector space $V=V_0 \oplus V_1$. The vectors of $V_0$ and $V_1$ are called homogeneous vectors of degree respectively equal to $0$ and $1$. For an homogeneous vector $x$, we denote by $\mid x\mid$ its degree.
A $\K$-Poisson superalgebra is a $\K$-super vector space  $\p=\p_0 \oplus \p_1$ equipped with two bilinear products denoted by $x\bu y$  and $\{x ,y \}$, having the following properties:
\begin{enumerate}
\item The couple $(\p,\bu )$ is a associative super commutative algebra,
that is, $$y \bu  x=(-1)^{\mid x \mid \mid y \mid}x \bu  y.$$ 
\item The couple $(\p, \{ , \})$ is  a Lie superalgebra, that is,
$$\{x ,y \}=-(-1)^{\mid x \mid \mid y \mid}\{y ,x \}$$
and satisfying the super Jacobi condition:
$$(-1)^{\mid z \mid \mid x \mid}\{x , \{y ,z \}\}+(-1)^{\mid x \mid \mid y \mid}\{y , \{z ,x \}\}+
(-1)^{\mid y \mid \mid z \mid}\{z , \{x ,y \}\}=0.$$
\item The products $\bu $ and $\{, \}$ satisfy the super Leibniz rule:
$$\{x,y\bu  z\}=\{x,y\} \bu  z+(-1)^{\mid x \mid \mid y \mid}y \bu  \{x,z\}.$$
where $x,y$ and $z$ are homogeneous vectors.
\end{enumerate}

\medskip
Let $(\p,\bu ,\{,\})$ be  a Poisson superalgebra.
Consider the multiplication
$$\begin{array}{l}
\mu(x,y)=xy=x \bu  y + \{ x,y \}.
\end{array}$$
We deduce 
$$
\left\{
\begin{array}{l}
\displaystyle x \bu  y=\frac{1}{2}(xy+(-1)^{\mid x \mid \mid y \mid}yx).\\
\displaystyle \{x , y\}=\frac{1}{2}(xy-(-1)^{\mid x \mid \mid y \mid}yx)
\end{array}
\right.
$$
The multiplication $\mu$ is called the depolarization of the couple $(\bu, \{,\})$.

The associativity condition writes for homogeneous vectors
$$\begin{array}{l}
\medskip
(x_1x_2)x_3+(-1)^{\mid x_1 \mid \mid x_2 \mid}(x_2x_1)x_3-(-1)^{\mid x_1 \mid \mid x_2 \mid+\mid x_1 \mid \mid x_3 \mid+\mid x_2 \mid \mid x_3 \mid}(x_3x_2)x_1\\ -(-1)^{\mid x_1 \mid \mid x_2 \mid+\mid x_1 \mid \mid x_3\mid} (x_2x_3)x_1\\
-x_1(x_2x_3)+(-1)^{\mid x_1 \mid \mid x_2 \mid+\mid x_1 \mid \mid x_3 \mid+\mid x_2 \mid \mid x_3 \mid}x_3(x_2x_1)-(-1)^{\mid x_2 \mid \mid x_3 \mid}x_1(x_3x_2)\\+(-1)^{\mid x_1 \mid \mid x_3 \mid+\mid x_2 \mid \mid x_3 \mid}x_3(x_1x_2)=0
\end{array}$$
Of course, we find again the relation of the non graded case considering all the elements of degree $0$. Using the same notations, we have 
$$W_1=^t(1, (-1)^{\mid x_1 \mid \mid x_2 \mid},-(-1)^{\mid x_1 \mid \mid x_2 \mid+\mid x_1 \mid \mid x_3 \mid+\mid x_2 \mid \mid x_3 \mid},0, -(-1)^{\mid x_1 \mid \mid x_2 \mid+\mid x_1 \mid \mid x_3 \mid} ,0)$$
$$W_2=^t(-1,0,(-1)^{\mid x_1 \mid \mid x_2 \mid+\mid x_1 \mid \mid x_3 \mid+\mid x_2 \mid \mid x_3 \mid},-(-1)^{\mid x_2 \mid \mid x_3 \mid},0,(-1)^{\mid x_1 \mid \mid x_3 \mid+\mid x_2 \mid \mid x_3 \mid})$$

The super Jacobi condition is written
$$\begin{array}{l}
\medskip
(-1)^{\mid x_1 \mid \mid x_3 \mid}(x_1x_2)x_3-(-1)^{\mid x_1 \mid \mid x_2 \mid+\mid x_1 \mid \mid x_3 \mid}(x_2x_1)x_3-(-1)^{\mid x_1 \mid \mid x_2 \mid+\mid x_2 \mid \mid x_3 \mid}(x_3x_2)x_1\\
 -(-1)^{\mid x_1 \mid \mid x_3 \mid+\mid x_2 \mid \mid x_3 \mid}(x_1(x_3x_2)+(-1)^{\mid x_1 \mid \mid x_2 \mid}(x_2x_3)x_1+(-1)^{\mid x_2 \mid \mid x_3 \mid}(x_3x_1)x_2 \\
-(-1)^{\mid x_1 \mid \mid x_3 \mid}x_1(x_2x_3)+(-1)^{\mid x_1 \mid \mid x_2 \mid+\mid x_1 \mid \mid x_3 \mid}x_2(x_1x_3)+(-1)^{\mid x_1 \mid \mid x_2 \mid+\mid x_2 \mid \mid x_3 \mid}x_3(x_2x_1)\\
+(-1)^{\mid x_1 \mid \mid x_3 \mid+\mid x_2 \mid \mid x_3 \mid}(x_1(x_3x_2)-(-1)^{\mid x_1 \mid \mid x_2 \mid}x_2(x_3x_1)-(-1)^{\mid x_2 \mid \mid x_3 \mid}x_3(x_1x_2 )=0
\end{array}$$
This relation is associated with the vectors
$$W_3=
\begin{pmatrix}
     (-1)^{\mid x_1 \mid \mid x_3 \mid}  \\
     -(-1)^{\mid x_1 \mid \mid x_2 \mid+\mid x_1 \mid \mid x_3 \mid}\\   
     -(-1)^{\mid x_1 \mid \mid x_2 \mid+\mid x_2 \mid \mid x_3 \mid}\\
      -(-1)^{\mid x_1 \mid \mid x_3 \mid+\mid x_2 \mid \mid x_3 \mid}\\
      (-1)^{\mid x_1 \mid \mid x_2 \mid}\\
 (-1)^{\mid x_2 \mid \mid x_3 \mid}
\end{pmatrix}, 
W_4=\begin{pmatrix}
    (-(-1)^{\mid x_1 \mid \mid x_3 \mid}\\
   (-1)^{\mid x_1 \mid \mid x_2 \mid+\mid x_1 \mid \mid x_3 \mid}\\
   (-1)^{\mid x_1 \mid \mid x_2 \mid+\mid x_2 \mid \mid x_3 \mid}\\
   (-1)^{\mid x_1 \mid \mid x_3 \mid+\mid x_2 \mid \mid x_3 \mid}\\
   -(-1)^{\mid x_1 \mid \mid x_2 \mid}\\
   -(-1)^{\mid x_2 \mid \mid x_3 \mid}
\end{pmatrix}
$$
From the study of the non graded case, we can consider that $B=Id$. Let  $A_g$ be the "graded" matrix $A$ associated with a graded relation satisfied by $\mu$.
For any triple of homogeneous vectors $x_1,x_2,x_3$, we denote by $(\alpha,\beta,\gamma)$ the element $(-1)^{\alpha\mid x_1 \mid+\beta \mid x_2 \mid+\gamma\mid x_3\mid}$. Any coefficient of a graded relation which defines $\m$ can be written $a_i(\alpha_i,\beta_i,\gamma_i)$ where $a_i$ will be the coefficient corresponding to the non graded case, or for $x_1,x_2,x_3$ homogeneous of degree $0$. Thus the matrix $A_g$ is
$$\begin{pmatrix}
    a_1(\alpha_1,\beta_1,\gamma_1)  &   a_2(\alpha_2,\beta_2,\gamma_2)  & a_3(\alpha_3,\beta_3,\gamma_3)  &  a_4(\alpha_4,\beta_4,\gamma_4)  &  a_6(\alpha_6,\beta_6,\gamma_6)  &  a_5(\alpha_5,\beta_5,\gamma_5)      \\
   a_2(\alpha_2,\gamma_2,\beta_2)    &  a_1(\alpha_1,\gamma_1,\beta_1)  &a_5(\alpha_5,\gamma_5,\beta_5) &a_6(\alpha_6,\gamma_6,\beta_6) &a_4(\alpha_4,\gamma_4,\beta_4)  &a_3(\alpha_3,\gamma_3,\beta_3) \\
   a_3(\gamma_3,\beta_3,\alpha_3)  &  a_6(\gamma_6,\beta_6,\alpha_6)  &a_1(\gamma_1,\beta_1,\alpha_1)  & a_5(\gamma_5,\beta_5,\alpha_5)   & a_2(\gamma_2,\beta_2,\alpha_2)  &a_4(\gamma_4,\beta_4,\alpha_4)\\
   a_4(\beta_4,\alpha_4,\gamma_4)&a_5(\beta_5,\alpha_5,\gamma_5)  &a_6(\beta_6,\alpha_6,\gamma_6)  &  a_1(\beta_1,\alpha_1,\gamma_1)  & a_3(\beta_3,\alpha_3,\gamma_3)  & a_2(\beta_2,\alpha_2,\gamma_2) \\
   a_5(\beta_5,\gamma_5,\alpha_5)  &a_4(\beta_4,\gamma_4,\alpha_4)& a_2(\beta_2,\gamma_2,\alpha_2)  & a_3(\beta_3,\gamma_3,\alpha_3)  & a_1(\beta_1,\gamma_1,\alpha_1)   &a_6(\beta_6,\gamma_6,\alpha_6) \\
   a_6(\gamma_6,\alpha_6,\beta_6) & a_3(\gamma_3,\alpha_3,\beta_3)  &a_4(\gamma_4,\alpha_4,\beta_4)& a_2(\gamma_2,\alpha_2,\beta_2)  &a_5(\gamma_5,\alpha_5,\beta_5)  &a_1(\gamma_1,\alpha_1,\beta_1) 
\end{pmatrix}
$$
where to shorten the writing a little more, we have l noted
$$ a_i(\alpha_i,\beta_i,\gamma_i) =a_i(-1)^{\alpha\mid x_i \mid+\beta \mid x_i \mid+\gamma\mid x_i\mid}$$
with $\alpha_i,\beta_i,\gamma_i \in \{0,1\}$
The associativity of $\bu$ is given by the solutions of the linear system
$$
\begin{array}{l}
   a_1(\alpha_1,\beta_1,\gamma_1) +(-1)^{123}a_3(\alpha_3,\beta_3,\gamma_3) -(-1)^{23}a_4(\alpha_4,\beta_4,\gamma_4)+(-1)^{13+23}a_5(\alpha_5,\beta_5,\gamma_5)=1  \\
    a_2(\alpha_2,\gamma_2,\beta_2)    +(-1)^{123}a_5(\alpha_5,\gamma_5,\beta_5) -(-1)^{23}a_6(\alpha_6,\gamma_6,\beta_6) +(-1)^{13+23}a_3(\alpha_3,\gamma_3,\beta_3) =(-1)^{12} \\
      a_3(\gamma_3,\beta_3,\alpha_3)  +(-1)^{123}a_1(\gamma_1,\beta_1,\alpha_1)  -(-1)^{23} a_5(\gamma_5,\beta_5,\alpha_5) +(-1)^{13+23}a_4(\gamma_4,\beta_4,\alpha_4)=-(-1)^{123} \\
         a_4(\beta_4,\alpha_4,\gamma_4)+(-1)^{123}a_6(\beta_6,\alpha_6,\gamma_6)   -(-1)^{23}a_1(\beta_1,\alpha_1,\gamma_1)  + (-1)^{13+23} a_2(\beta_2,\alpha_2,\gamma_2) =0 \\
   a_5(\beta_5,\gamma_5,\alpha_5)  +(-1)^{123} a_2(\beta_2,\gamma_2,\alpha_2)  -(-1)^{23} a_3(\beta_3,\gamma_3,\alpha_3)  +(-1)^{13+23}a_6(\beta_6,\gamma_6,\alpha_6) =-(-1)^{12+13} \\
                a_6(\gamma_6,\alpha_6,\beta_6) +(-1)^{123}a_4(\gamma_4,\alpha_4,\beta_4) -(-1)^{23} a_2(\gamma_2,\alpha_2,\beta_2) +(-1)^{13+23}a_1(\gamma_1,\alpha_1,\beta_1) =0 \\
   \end{array}
$$
where 
$$(-1)^{123}=(-1)^{\mid x_1 \mid \mid x_2 \mid+\mid x_1 \mid \mid x_3 \mid+\mid x_2 \mid \mid x_3 \mid}, (-1)^{23}=(-1)^{\mid x_2 \mid \mid x_3 \mid}, (-1)^{13+23}=(-1)^{\mid x_1 \mid \mid x_3 \mid+\mid x_2 \mid \mid x_3 \mid}.$$
If the homogeneous vectors $x_i$ are of degree $0$, then we find again the system defined in the non graded case which has the solution
$$a_6=a_1-a_2+a_4, \ a_5=1+a_1-a_3+a_4.$$
Similarly Jacobi’s conditions in degree $0$ lead to
$$a_4=-2-3a_1+2a_2+2a_3.$$
We deduce that $\mu$ satisfies a relation of type
$$
\begin{array}{l}
 \widetilde{a_1}(x_1x_2)x_3+\widetilde{a_2 }(x_2x_1)x_3+\widetilde{a_3}(x_3x_2)x_1+(-2-3\widetilde{a_1}+2\widetilde{a_2}+2\widetilde{a_3})(x_1x_3)x_2\\
+(-1-2\widetilde{a_1}+2\widetilde{a_2}+\widetilde{a_3})(x_2x_3)x_1+(-2-2\widetilde{a_1}+\widetilde{a_2}+2\widetilde{a_3})(x_3x_1)x_2-x_1(x_2x_3)=0
\end{array}
$$
where $\widetilde{a_i}=a_i(\alpha_i,\beta_i,\gamma_i).$
Now the study of the graded -Leibniz identity means to 
$$a_1=1,a_2=1/3,a_3=0,a_4=1/3,a_5=-1/3,a_6=-1/3.$$
When we look these equations when the degree of $x_1,x_2,x_3$ are respectively equal to $(1,1,0), $ $(1,0,1),(0,1,1),(1,1,1)$, we obtain the following result:

\begin{theorem}
Let $\p$ a $\K$-supervector space. Thus $\p$ is a Poisson superalgebra if and only if there exists on $\p$ a nonassociative product $x  y$ satisfying
 \begin{equation}\label{superPoisson}
\left\{
 \begin{array}{l}
\medskip
 3(xy)z-3x(yz)  + (-1)^{\mid x \mid \mid y \mid}(yx)z -(-1)^{\mid y \mid \mid z \mid} (xz)y
 - (-1)^{\mid x \mid \mid y \mid+\mid x \mid \mid z \mid}(yz)x \\
 + (-1)^{\mid x \mid \mid z \mid+\mid y \mid \mid z \mid}(zx)y =0 
\end{array}
\right.
\end{equation}
for any homogeneous vectors $x,y,z \in\p$.
\end{theorem}

\medskip
\noindent{\bf Examples.} Any $2$-dimensional superalgebra $\p=V_0 \oplus V_1$ with an homogeneous basis $\{e_0,e_1\}$ writes
$$\left\{
\begin{array}{l}
e_0e_0=ae_0,\\
 e_0e_1=be_1, \ e_1e_0=ce_1,\\
e_1e_1=de_0.
\end{array}
\right.
$$
This product is a superPoisson product if and only if we have
$$\left\{
\begin{array}{l}
d=0,\\
3(a-b)b+ab-2bc+c^2=0,\\
3(a-c)c+ab-2bc+c^2=0,
\end{array}
\right.
$$
or
$$\left\{
\begin{array}{l}
a=0,\\
b=-c.
\end{array}
\right.
$$
We obtain the following $2$-dimensional Poisson superalgebras
$$\left\{
\begin{array}{ll}
\mathcal{SP}_{2,1} & e_0e_0=ae_1\\
\mathcal{SP}_{2,2} & e_0e_0=ae_1, e_0e_1=e_1e_0=ae_1\\
\mathcal{SP}_{2,3} &  e_0e_1=-e_1e_0=be_1\\
\mathcal{SP}_{2,4} &  e_0e_1=-e_1e_0=be_1, \ e_1e_1=de_0,\\
\end{array}
\right.
$$
the non written product being considered equal to $0$.

\subsection{Properties of Poisson superalgebras }
\begin{definition}
A nonassociative superalgebra  is called superflexive if  the multiplication $xy$ satisfy
$$A(x,y,z) + (-1)^{(|x||z|+|x||y|+|y||z|)}A(z,y,x)=0$$
for any homogeneous elemnts $x,y,z$, where $A(x,y,z)=(xy)z-x(yz)$ is the associator of the multiplication.
\end{definition}
\begin{proposition}
Let $\p$ be a Poisson superalgebra. Then the nonassociative product defining the superPoisson structure is superflexive.
\end{proposition}
\noindent{\it Proof.} In fact, let 
$$B(x,y,z)=3(A(x,y,z) + (-1)^{(|x||z|+|x||y|+|y||z|)}A(z,y,x)) .
$$
We  have
$$
\begin{array}
{rl}
\medskip
B(x,y,z) =& -(-1)^{\mid x \mid \mid y \mid}(yx)z +(-1)^{\mid y \mid \mid z \mid} (xz)y
 + (-1)^{\mid x \mid \mid y \mid+\mid x \mid \mid z \mid}(yz)x 
 - (-1)^{\mid x \mid \mid z \mid+\mid y \mid \mid z \mid}(zx)y \\
 \medskip
 
 &+ (-1)^{(|x||z|+|x||y|+|y||z|)}(-(-1)^{\mid z \mid \mid y \mid}(yz)x +(-1)^{\mid y \mid \mid x \mid} (zx)y \\
 \medskip
 
 &+ (-1)^{\mid z \mid \mid y \mid+\mid z \mid \mid x \mid}(yx)z 
 - (-1)^{\mid z \mid \mid x \mid+\mid y \mid \mid x \mid}(xz)y )\\
 \medskip
 
 =& (-(-1)^{\mid x \mid \mid y \mid} +(-1)^{\mid x \mid \mid y \mid}(yx)z + ((-1)^{\mid y \mid \mid z \mid} -(-1)^{\mid y \mid \mid z \mid} (xz)y \\
 \medskip
 &+((-1)^{\mid x \mid \mid y \mid+\mid x \mid \mid z \mid} - (-1)^{\mid x \mid \mid y \mid+\mid x \mid \mid z \mid}(yz)x 
 +( - (-1)^{\mid x \mid \mid z \mid+\mid y \mid \mid z \mid}\\
 \medskip
 & +(-1)^{\mid x \mid \mid z \mid+\mid y \mid \mid z \mid})(zx)y \\
 \medskip
 =& 0.
 \end{array}
 $$
 
 \medskip
 
 \noindent{\bf Remark : on the power associativity.} Recall that a non-associative algebra is power associative if every element
generates an associative subalgebra. Let $\p$ be a Poisson superalgebra provided with its non associative product $xy$. If $V_0$ is its the even homogeneous part, then the 
restriction of the product $xy$ is a multiplication in 
this homogeneous vector space satisfying Identity \ref{superPoisson}. Since all the vectors of $V_0$ are of degree $0$, Identity \ref{superPoisson} is reduced to identity
 \ref{associator}. We deduce that $V_0$ is a Poisson algebra and any vector $x$ in $V_0$ generates an associative subalgebra of $V_0$ and of $\p$.
  
  Assume now that $y$ is an odd vector. We have
  $$y.y=\displaystyle\frac{1}{2}(yy+(-1)yy)=0,$$
  and
  $$\{y,y\}=\displaystyle\frac{1}{2}(yy-(-1)yy)=yy.$$
  If we write $y^2=yy$, then $$y^2=\{y,y\}.$$
  This implies
  $$yy^2=y\{y,y\}=y \bu  \{y,y\}+\{y,\{y,y\}\}.$$
  But from the suoer identity of Jacobi, $\{y,\{y,y\}\}=0.$ Thus we have
  $$yy^2= y \bu  \{y,y\}=\{y,y\} \bu  y=y^2y.$$
  We can write
  $$y^3=yy^2=y^2y.$$
  Now
  $$y^2y^2=\{y,y\}\{y,y\}=\{y,y\}\bu \{y,y\}+\{\{y,y\},\{y,y\}\}.$$
  We have also
  $$yy^3=y\bu  y^3+\{y,y^3\}=y\bu  y\bu \{y,y\}+\{y,y\bu  \{y,y\}\}.$$
  But $y \bu  y=0$. Thus, from the Leibniz rule,
  $$yy^3=\{y,y\bu  \{y,y\}\}=-y\bu  \{y,\{y,y\}\}+\{y,y\}\bu \{y,y\}=\{y,y\}\bu  \{y,y\}.$$
  We deduce
  $$y^2y^2-yy^3= \{\{y,y\},\{y,y\}\}.$$
  Since $\{y,y\}$ is of degree $0$, we obtain
$$y^2y^2-yy^3=0.$$
We can write
$$y^4=y^2y^2=yy^3=y^3y$$
the last equality resutlts of $\{y,y\bu  \{y,y\}\}=\{y\bu  \{y,y\},y\}.$ 
Now, using the Identity \ref{superPoisson} to the triple $(y^i, y^j,y^k)$ with $i+j+k=5$, we obtain a 
linear system 
on the vectors $y^iy^j$ with $i+j=5$, which admits as solutions
$$yy^4=y^2y^3= y^3y^2=y^4y.$$
Thus $y^5$ is well determinated. By induction, using Identity \ref{superPoisson}  on the triple $(y^i, y^j,y^k)$ with $i+j+k=n$, using induction hypothesis $y^py^{n-1-p}=y^{n-1}$, we obtain that
$$y^n=y^py^{n-p}$$
for any $p=1,\bu s,n-1$. Thus any homogeneous element of odd degree generates an associative algebra. But this property is not true, in general, for non homogeneous element. For example
if $x=e_0+e_1$ with degree of $e_i$ equal to $i$, then
$$\left\{
\begin{array}{lll}
x^2x-xx^2&=&A(e_1,e_1,e_0)+A(e_1,e_0,e_1)+A(e_0,e_1,e_1)\\
&&+A(e_1,e_0,e_0)+A(e_0,e_0,e_1)+A(e_0,e_1,e_0)\\
&=& 4(e_1^2e_0-(e_0e_1)e_1)
\end{array}
\right.
$$
which is not trivial.

\medskip
\subsection{Transposed Poisson superalgebras}
Let  $V=V_0 \oplus V_1$ be a $\Z_2$-graded vector space . T
A $\K$-transposed Poisson superalgebra is a $\K$-super vector space  $\p=\p_0 \oplus \p_1$ equipped with two bilinear products denoted by $x\bu y$  and $\{x ,y \}$, having the following properties:
\begin{enumerate}
\item The couple $(\p,\bu )$ is a associative supercommutative algebra,
that is, $$y \bu  x=(-1)^{\mid x \mid \mid y \mid}x \bu  y.$$ 
\item The couple $(\p, \{ , \})$ is  a Lie superalgebra, that is,
$$\{x ,y \}=-(-1)^{\mid x \mid \mid y \mid}\{y ,x \}$$
and satisfying the super Jacobi condition:
$$(-1)^{\mid z \mid \mid x \mid}\{x , \{y ,z \}\}+(-1)^{\mid x \mid \mid y \mid}\{y , \{z ,x \}\}+
(-1)^{\mid y \mid \mid z \mid}\{z , \{x ,y \}\}=0.$$
\item The products $\bu $ and $\{, \}$ satisfy the super transposed Leibniz rule:
 $$2x \bu [y, z] = [x \bu y, z] +(-1)^{\mid x \mid \mid y\mid} [y, x \bu z].$$
where $x,y$ and $z$ are homogeneous vectors.
\end{enumerate}

\medskip
Let $(\p,\bu ,\{,\})$ be  a transposed Poisson superalgebra.
Consider the multiplication
$$\begin{array}{l}
\mu(x,y)=xy=x \bu  y + \{ x,y \}.
\end{array}$$
We deduce 
$$
\left\{
\begin{array}{l}
\displaystyle x \bu  y=\frac{1}{2}(xy+(-1)^{\mid x \mid \mid y \mid}yx).\\
\displaystyle \{x , y\}=\frac{1}{2}(xy-(-1)^{\mid x \mid \mid y \mid}yx)
\end{array}
\right.
$$
The multiplication $\mu$ is called the depolarization of the couple $(\bu, \{,\})$. 
\begin{proposition}
The depolarization $\mu$ of the couple $(\bu,[-,-]$ which defines a structure of Transposed Poisson superalgebra on $V$ have to satisfy the two independent axioms
\begin{enumerate}
  \item $$
\begin{array}{l}
a_1(xy)z+(-1)^{\mid x \mid \mid y \mid}a_2(yx)z+(-1)^{\mid x \mid \mid z \mid+\mid y \mid \mid z \mid}a_3(zy)x+(-1)^{\mid y \mid \mid z\mid}(-2-3a_1+2a_2+2a_3)(xz)y\\
+(-1)^{\mid x \mid \mid y \mid +\mid x \mid \mid z \mid}(-1-2a_1+2a_2+a_3)(yz)x\\+(-1)^{\mid y \mid \mid z \mid+\mid x \mid \mid z \mid}((-2-2a_1+a_2+2a_3)(zx)y+x(yz)=0.
\end{array}
$$
\item
$$\begin{array}{l}
      2(xy)z  -x(yz)-(-1)^{\mid x \mid \mid y \mid}(2(yx)z-y(xz))+(-1)^{\mid x \mid \mid z \mid+\mid y \mid \mid z \mid}((zy)x-2z(yx))\\-(-1)^{\mid y \mid \mid z\mid}(xz)y+(-1)^{\mid x \mid \mid y \mid +\mid x \mid \mid z \mid}(yz)x
      -(-1)^{\mid y \mid \mid z \mid+\mid x \mid \mid z \mid}(zx)y  \\
   -x(zy)+y(zx)+2z(xy)=0  
\end{array}
$$
where $(x,y,z)$ are hogeneous vectors.
\end{enumerate}
\end{proposition}

\section{Other depolarizations}

In \cite{HLS} one introduce a new type of algebras defining by two multiplications and a distributive law, for expand some deformations and called Hom-Lie algebras. Such algebra $V$ is defined by an anticommutative multiplication $[-,-]$, an endomorphism $f$ of $V$ satisfying
$$[[x_1,x_2],f(x_3)]+[[x_2,x_3],f(x_1)]+[[x_3,x_1],f(x_2)]=0.$$
Let us consider the multiplication $\bu$ given by 
$$x \bu y=[x,f(y)].$$
This multiplication is commutative if $f$ satisfies
$$[x,f(y)]=[y,f(x)]$$
that is, since the bracket is anticommutative
\begin{equation}
\label{gV}
[f(x),y]+[x,f(y)]=0
\end{equation}
for any $x,y \in V$.
Let $\mathcal{G}(V)$ the linear subspace of the Lie algebra $\mathcal{E}nd(V)$ whose elements satisfy (\ref{gV}). This space contains all the derivations of $(V,[-,-])$ which are trivial on the derived subalgebra. 
\begin{lemma}
$\mathcal{G}(V)$ is a Lie subalgebra of $\mathcal{E}nd(V)$.
\end{lemma}
\pf Let $f,b \in \mathcal{G}(V)$. Then
$$\begin{array}{ll}
 [[f,g](x),y]+[x,[f,g](y]     & =  [f(g(x),y]-[g(f(x),y]+[x,f(g(y))]-[x,g(f(y)] \\
      &   =[f(g(x),y]+[x,f(g(y))]-[g(f(x),y]-[x,g(f(y)] \\
      =0
\end{array}
$$
\noindent{\bf Example} If  $(V,[-,-])$ is the Heisenberg Lie algebra $\frak{h}_3$ of dimension $3$, then in the basis $\{e_1,e_2,e_3\}$ which satisfies $[e_1,e_2]=0$, then the element of $\mathcal{G}(\frak{h}_3)$ are given by the matrices
$$\begin{pmatrix}
     a_1 & a_2 & 0   \\
      b_1& -a_1 & 0\\
      c_1 & c_2 & c_3 
\end{pmatrix}
$$

\begin{lemma}
The commutative multiplication $\bu$ satisfies the following identities
$$\mathcal{A}_\bu (x,y,z)+\mathcal{A}_\bu (y,z,x)+\mathcal{A}_\bu (z,x,y)=0$$
$$\mathcal{A}_\bu (x,y,z)+\mathcal{A}_\bu (z,y,x)=0$$
that is the algebra $(V,\bu)$ is a commutative $(Id+c+c^2)$-algebra and a $(Id+\tau_{13})$-algebra (\cite{GRNonass}).
\end{lemma}
\pf In fact these identities are always true as soon as the multiplication is commutative.
\medskip

If $f \in \mathcal{G}(V)$, then on $V$ we have

\begin{enumerate}
  \item An anticommutative multiplication $[-,-]$,
  \item A commutative multiplication $\bu$,
  \item A distributive law
  $$[x_1,x_2] \bu x_3 +[x_2,x_3] \bu x_1 + [x_3,x_1] \bu x_2 =0.$$
\end{enumerate}
Let $\mu$ the multiplication
$$xy=\frac{1}{2}(x\bu y +[x,y]).$$

\begin{theorem}
The multiplication $\mu$ is a depolarization of the pair $(\bu,[-,-])$ if it satisfies the following relation
$$\begin{array}{l}\mathcal{AA}(x_1,x_2,x_3)-\mathcal{AA}(x_2,x_1,x_3)-\mathcal{AA}(x_3,x_2,x_1)-\mathcal{AA}(x_1,x_3,x_2)\\
+\mathcal{AA}(x_2,x_3,x_1)+\mathcal{AA}(x_3,x_1,x_2)=0
\end{array}$$
where $\mathcal{AA}(x_1,x_2,x_3)=(x_1x_2)x_3+x_1(x_2x_3)$ is the antiassociator of the multiplication $\mu$.
\end{theorem}

\end{document}